\documentclass[12pt,reqno]{amsart}

\usepackage{amscd}
\usepackage{amssymb}

\newcommand{\N}{{\mathbb N}}
\newcommand{\Z}{{\mathbb Z}}

\newcommand{\C}{{\mathbb C}}

\renewcommand{\P}{{\mathbb P}}

\newcommand{\AAA}{{\mathcal A}}

\newcommand{\DD}{{\mathcal D}}
\newcommand{\EE}{{\mathcal E}}

\newcommand{\II}{{\mathcal I}}
\newcommand{\KK}{{\mathcal K}}

\newcommand{\NN}{{\mathcal N}}

\newcommand{\OO}{{\mathcal O}}

\newcommand{\RR}{{\mathcal R}}
\newcommand{\TT}{{\mathcal T}}
\newcommand{\UU}{{\mathcal U}}
\newcommand{\VV}{{\mathcal V}}
\newcommand{\XX}{{\mathcal X}}

\newcommand{\ddd}{{\rm d}}
\newcommand{\www}{\widetilde}

\newcommand{\uuuu}{\underline}

\newcommand{\paa}{\partial}
\newcommand{\nnn}{\nabla}

\DeclareMathOperator{\Aut}{Aut}

\DeclareMathOperator{\id}{id}

\DeclareMathOperator{\Lie}{Lie}

\DeclareMathOperator{\rk}{rk}

\DeclareMathOperator{\tr}{tr}

\theoremstyle{plain}
\newtheorem{lemma}{Lemma}[section]
\newtheorem{theorem}[lemma]{Theorem}

\newtheorem{corollary}[lemma]{Corollary}
\newtheorem{conjecture}[lemma]{Conjecture}

\theoremstyle{definition}
\newtheorem{definition}[lemma]{Definition}

\newtheorem{remark}[lemma]{Remark}
\newtheorem{remarks}[lemma]{Remarks}
\newtheorem{example}[lemma]{Example}

\newtheorem{definition/lemma}[lemma]{Definition/Lemma}

\begin{document}

\title[Meromorphic connections over $F$-manifolds]
{Meromorphic connections over $F$-manifolds} 

\author[L. David and C. Hertling]{Liana David and Claus Hertling}
\address{Institute of Mathematics Simion Stoilow
of the Romanian Academy, Calea
Grivitei no. 21, Sector 1, Bucharest, Romania}
\email{liana.david@imar.ro}
\address{Lehrstuhl f\"ur Mathematik VI,
Universit\"at Mannheim, A5,6, 68131 Mannheim, Germany}
\email{hertling\char64 math.uni-mannheim.de}

\dedicatory{In remembrance of Boris Dubrovin}

\thanks{This work was funded by the Deutsche Forschungsgemeinschaft
(DFG, German Research Foundation) -- 242588615}



\date{December 04, 2019}

\begin{abstract}
This paper review one construction of Frobenius manifolds
(and slightly weaker structures). It splits it into several steps
and discusses the freedom and the constraints in these steps. 
The steps pass through holomorphic bundles with meromorphic connections.
A conjecture on existence and uniqueness of certain such bundles,
a proof of the conjecture in the 2-dimensional cases,  
and some other new results form a research part of this paper.
\end{abstract}

\maketitle


\tableofcontents


\section{Introduction}\label{c1}
\setcounter{equation}{0}

\noindent
We owe to Boris Dubrovin the notion of a {\it Frobenius manifold} \cite{Du92}. 
Frobenius manifolds are rich geometric objects, which have many different
facets and which are at the crossroads of very different mathematical areas.
To describe them, the language of differential geometry is good.
To understand them well, one has to consider also holomorphic vector bundles
with meromorphic connections.
They arise in integrable system, in quantum cohomology and in singularity
theory, there in the study of families of holomorphic functions.
Isomorphy of certain Frobenius manifolds from quantum cohomology and
certain Frobenius manifolds from singularity theory is one version of
mirror symmetry.

This paper review one construction of Frobenius manifolds
(and slightly weaker structures). It splits it into several steps
and discusses the freedom and the constraints in these steps. 
The steps pass through holomorphic bundles with meromorphic connections.
A conjecture on existence and uniqueness of certain such bundles,
a proof of the conjecture in the 2-dimensional cases,  
and some other new results form a research part of this paper.

A differential geometric description of a {\it Frobenius manifold
with Euler field} $(M,\circ,e,g,E)$ is as follows (Definition \ref{t3.1} (d)).
$M$ is a complex manifold $M$ with a (holomorphic) commutative and 
associative multiplication $\circ$ on the holomorphic tangent bundle $TM$.
$e$ is a unit field, a global holomorphic vector field with $e\circ=\id$.
$g$ is a {\it metric}, which means here a (holomorphic) bilinear
symmetric and nondegenerate pairing with values in $\C$ on $TM$.
Its Levi-Civita connection $D$ must be flat, the unit field must be $D$-flat,
$D(e)=0$, the metric must be multiplication invariant,
\begin{eqnarray}\label{1.1}
g(X\circ Y,Z)&=& g(X,Y\circ Z)\quad\textup{for }X,Y,Z\in\TT_M
\end{eqnarray}
($\TT_M$ is the sheaf of holomorphic vector fields),
and the connection $D$ must satisfy the potentiality condition
\begin{eqnarray}\label{1.2}
D_X(Y\circ) = D_Y(X\circ)&=& [X,Y]\circ \quad\textup{for }X,Y,Z\in\TT_M.
\end{eqnarray}
Finally, the Euler field is a global holomorphic vector field $E$
with 
\begin{eqnarray}\label{1.3}
\Lie_E(\circ)=1\cdot \circ ,\quad \Lie_E(g)=(2-d)\cdot g
\quad\textup{for some }d\in\C.
\end{eqnarray}

\eqref{1.2} is called {\it potentiality condition}, because
together with the flatness of $D$ and \eqref{1.1} it implies
locally the existence of a potential $F\in\OO_M$ with
\begin{eqnarray}\label{1.4}
XYZ(F)=g(X\circ Y,Z)\quad\textup{for }D\textup{-flat vector fields }X,Y,Z.
\end{eqnarray}
The associativity of the multiplication $\circ$ is then expressed
in the WDVV-equations for $F$. 
This potential is the central object in the construction of 
Frobenius in quantum cohomology. 

But in the construction in singularity theory, it plays no role.
There one starts with a manifold with multiplication and unit field
and Euler field $(M,\circ,e,E)$ and then {\it adds} a metric $g$.
This {\it adding a metric} is highly nontrivial. 
One aim of this paper is to review this procedure, to split it up
into several steps and to discuss the freedom and the constraints 
in these steps. They go through other holomorphic bundles on $M$
(not the tangent bundle $TM$) with meromorphic connections.

The manifold with multiplication $(M,\circ,e)$ inherits from the
potentiality condition a specific integrability condition of the
multiplication, see \eqref{2.1}. With this condition, it is 
called an {\it F-manifold} \cite{HM99}. An Euler field on an
F-manifold is a global holomorphic vector field $E$ with
$\Lie_E(\circ)=\circ$. 

We will not only consider Frobenius manifolds with Euler fields,
but also slightly weaker structures, altogether 4 cases:
Frobenius manifolds without/with Euler fields and without/with metrics.
Frobenius manifolds without Euler field and without metric
(but with flat and torsion free connection $D$ with potentiality)
were considered in \cite{Ma05} and called {\it flat F-manifolds}. 
So, we consider the four cases, 

\medskip
\begin{tabular}{c|c|c}
 & without Euler field & with Euler field \\ \hline
without metric & (a) flat F-manifold & (b) flat F-manifold with $E$\\
\hline
with metric & (c) Frobenius mfd. & (d) Frobenius mfd. with $E$
\end{tabular}

\medskip
New interest especially in the flat F-manifolds with Euler fields
arose in recent work.
Dubrovin \cite{Du96} had observed a relation between
3-dimensional semisimple Frobenius manifolds with Euler fields
and the solutions of a 1-parameter subfamily of the Painlev\'e
VI equations. Beautiful generalizations were found in 
\cite{Lo14} \cite{KMS15} \cite{AL15} \cite{KM19}, relations between
all generic 3- and 4-dimensional {\it regular} F-manifolds
and all generic solutions of the Painlev\'e equations of types
VI, V, IV, III and II (see the Remarks \ref{t3.5}).
Also, the constructions of Frobenius manifolds with Euler fields
on the orbit spaces of the Coxeter groups \cite{SaK79} \cite{Du96} 
were generalized to constructions of flat F-manifolds with
Euler fields on the orbit spaces of most finite complex
reflection groups \cite{KMS15} \cite{KMS18} \cite{AL17}
(see the Remarks \ref{t3.4}).

The construction which we want to discuss proceeds in the case
of a flat F-manifold with Euler field as follows.
One has to build up structures in 4 steps:

\begin{list}{}{}
\item[(I)(b)] An F-manifold $M$ with Euler field.
\item[(II)(b)] A $(TE)$-structure over the F-manifold $M$ with
Euler field.
\item[(III)(b)] An extension of the $(TE)$-structure to a 
pure $(TLE)$-structure.
\item[(IV)(b)] A choice of a primitive section which allows
to shift a part of the connection of the $(TE)$-structure
to the F-manifold. The resulting flat connection enriches
the F-manifold to a flat F-manifold.
\end{list}

This recipe is discussed in detail in the Remarks \ref{t6.11}.
The ingredients are treated before, F-manifolds with Euler 
fields in section \ref{c2}, $(TE)$-structures in section 
\ref{c4}, pure $(TLE)$-structures in section \ref{c5},
and primitive sections in \ref{c6}.

A variant of this recipe was used first by Morihiko Saito 
in \cite{SaM89} in the case of isolated hypersurface singularities.
The base space of a universal unfolding is an F-manifold
with Euler field \cite{He02}. The universal unfolding comes equipped 
with a vector bundle with a meromorphic connection with
logarithmic poles along a discriminant, 
the Gauss-Manin connection (e.g. \cite{He02}).
A partial Fourier-Laplace transformation gives a
$(TE)$-structure. \cite{SaM89} starts essentially from this 
situation and shows that the steps (III)(b) and (IV)(b)
can be carried out. Here (III)(b) amounts to an extension
of a bundle on $\C$ with a meromorphic connection with a 
pole at 0 to a trivial bundle on $\P^1$ with a logarithmic
pole at $\infty$. This is a Birkhoff problem.
The solution in \cite{SaM89} starts with a Hodge filtration
which is induced from the pole at 0 and chooses a
monodromy invariant opposite filtration (which indeed exists).
That filtration induces the desired extension to $\infty$.
(IV)(b) amounts to the existence of a global section with
certain properties. 

In fact, \cite{SaM89} works with the Gauss-Manin connection,
not with $(TE)$-structures. His construction is also
covered in \cite{He02}. 
The recipe as above (and a lot of background material)
is discussed in Claude Sabbah's book \cite{Sa02}.

In this paper we do not treat the Gauss-Manin connection
and its partial Fourier-Laplace transform, but we take
an abstract point of view and study $(TE)$-structures in 
general. 

The freedom and the constraints in the 4 steps are
surprisingly different in the 4 cases (a), (b), (c) and (d).
The Remarks \ref{t6.11} discuss this in detail.
The following diagram comprises the situation in a facile way
in the cases (a) and (b). It shows in the left column several steps
from an F-manifold (without Euler field) to a flat F-manifold
(without Euler field), and in the right column the same steps
from an F-manifold with Euler field to a flat F-manifold with Euler field.

\vspace*{0.4cm}
\noindent
\begin{tabular}{ll}
\fbox{\parbox{5cm}{(I)(a) F-manifold: $\exists$ easy}}  & 
\fbox{\parbox{6.5cm}{(I)(b) F-manifold with Euler field: $\exists$ easy}} 
\end{tabular}

\noindent
\begin{tabular}{ll}
\fbox{\parbox{5cm}{(II)(a) (T)-structure over an F-manifold:\\
$\exists$ in the gen. semisimple cases: Conjecture \ref{t7.2} (a)\\
$\exists$ in general: no \\
If $\exists$, possibly functional par.}}&
\fbox{\parbox{6.5cm}{(II)(b) (TE)-structure over an F-manifold
with Euler field: \\
$\exists$ in the generically semisimple cases: Conjecture \ref{t7.2} (b)\\
$\exists$ in general: no\\ 
If $\exists$, functional par.}}
\end{tabular}

\noindent
\begin{tabular}{ll}
\fbox{\parbox{5cm}{(III)(a) pure (TL)-structure over an F-manifold:\\
$\exists$ for free, functional par.:\\
Theorem \ref{t5.1} (a)}} & 
\fbox{\parbox{6.5cm}{(III)(b) pure (TLE)-structure over an F-manifold
with Euler field:
$\exists$ is a Riemann-Hilbert-Birkhoff problem}} 
\end{tabular}

\noindent
\begin{tabular}{ll}
\fbox{\parbox{5cm}{(IV)(a) flat F-manifold:\\
$\exists$ for free, $<\infty$ par.}} & 
\fbox{\parbox{6.5cm}{(IV)(b) flat F-manifold with Euler field:
$\exists$ needs $\exists$ of a primitive section which is
an eigenvector of $Q$, often
yes, but not always, $<\infty$ par.}}
\end{tabular}

\vspace*{0.4cm}
Section \ref{c7} proposes Conjecture \ref{t7.2}. 
It makes a precise prediction about existence and uniqueness 
and the parameters for all $(TE)$-structures over an
irreducible germ of a generically semisimple F-manifold. 
Section \ref{c8} proves this conjecture in the 
2-dimensional cases. There one has the F-manifolds $I_2(m)$
($m\in\Z_{\geq 3}$). The proof requires quite some diligence
and uses a formal classification of $(T)$-structures
over $I_2(m)$ in \cite{DH19-1}.

The paper is structured as follows.
Section \ref{c2} discusses old and new results on F-manifolds,
with and without Euler fields.
Section \ref{c3} introduces Frobenius manifolds without/with Euler
fields and without/with metrics and reports on the recent
interest in flat F-manifolds with Euler fields. 
Section \ref{c4} defines $(TE)$-structures and related structures.
It serves just as a dictionary.
Section \ref{c5} treats the problem of extending $(TE)$-structures
to pure $(TLE)$-structures, a Birkhoff problem, and it states 
the possibly surprising fact that the analogous extension problem
for $(T)$-structures has always many solutions.
Section \ref{c6} gives the recipe in 4 steps and discusses them.
Section \ref{c7} proposes Conjecture \ref{t7.2}.
Section \ref{c8} solves it for the F-manifolds of types $I_2(m)$.

This paper is a mixture of a survey and a research paper. 
It is related to the book \cite{Sa02}, but it emphasizes different points.
Conjecture \ref{t7.2} and section \ref{c8} are new. 
Some new results are also in the sections \ref{c5} and \ref{c6}
(Theorem \ref{t5.1} (a)+(b), Theorem \ref{t5.6}, 
Corollary \ref{t5.7}, Theorem \ref{t6.7}). 
Other classical references on Frobenius manifolds
are \cite{Du96}, \cite{Ma99} and \cite{He02}.

Some notation:
$\N=\{1,2,...\}$, $\N_0=\N\cup\{0\}$.
$M$ is always a complex manifold,
$TM$ is its holomorphic tangent bundle, and 
$\TT_M$ is the sheaf of holomorphic sections of $TM$,
so the sheaf of holomorphic vector fields.
$\circ$ denotes always a holomorphic commutative
and associative multiplication on $TM$.
And then $e\in\TT_M$ is a global holomorphic vector field
with $e\circ=\id$.
Usually local coordinates on $M$ will be called
$t=(t_1,...,t_n)$, and their coordinate vector fields
will be written as $\paa_k=\paa/\paa t_k$.

\section{A review of F-manifolds}\label{c2}
\setcounter{equation}{0}

\noindent 
$F$-manifolds were first defined in \cite{HM99}.
This section reviews their definition, basic properties
which were developed in \cite{He02},
and some younger results.

\begin{definition}\label{t2.1}
(a) An {\it F-manifold} $(M,\circ,e)$ (without Euler field) 
is a holomorphic manifold $M$ with a holomorphic 
commutative and associative multiplication $\circ$
on the holomorphic tangent bundle $TM$ and with a 
global holomorphic vector field $e\in\TT_M$ with
$e\circ=\id$ ($e$ is called a {\it unit field}),
which satisfies the following integrability condition:
\begin{eqnarray}\label{2.1}
\Lie_{X\circ Y}(\circ)&=& X\circ\Lie_Y(\circ)+Y\circ\Lie_X(\circ)
\textup{ for }X,Y\in\TT_M.
\end{eqnarray}
(b) Given an F-manifold $(M,\circ,e)$, an {\it Euler field} on it is a global
vector field $E\in\TT_M$ with $\Lie_E(\circ)=\circ$.
\end{definition}

At first sight, the integrability condition \eqref{2.1}
might look surprising. Below in Theorem \ref{t2.12}, 
an equivalent condition is given. 
The integrability condition is also important in the proof 
in \cite{He02} of the basic decomposition result
Theorem \ref{t2.3} on germs of F-manifolds.
There it allows to extend the pointwise decomposition
of $T_tM$ in Remark \ref{t2.2} into a local decomposition
of the F-manifold.

\begin{remark}\label{t2.2}
A finite dimensional commutative and associative $\C$-algebra
$A$ with unit $e\in A$ decomposes uniquely into a direct sum
$A=\bigoplus_{k=1}^l A_k$ of local and irreducible
algebras $A_k$ with units $e_k$ with $e=\sum_{k=1}^l e_k$
and $A_{k_1}\circ A_{k_2}=0$ for $k_1\neq k_2$.
This is elementary (linear) algebra. The
decomposition is obtained as the simultaneous decomposition
into generalized eigenspaces of all endomorphisms
$a\circ:A\to A$ for $a\in A$ 
(see e.g. Lemma 2.1 in \cite{He02}).
The algebra $A$ is called {\it semisimple} if
$l=\dim A$ (so then $A_k=\C\cdot e_k$ for all $k$).
\end{remark}

\begin{theorem}\label{t2.3}\cite[Theorem 2.11]{He02}
Let $((M,t^0),\circ,e)$ be the germ at $t^0$ of an
F-manifold. 

\medskip
(a) The decomposition of the algebra 
$(T_{t^0}M,\circ|_{t^0},e|_{t^0})$ with unit into local
algebras extends into a canonical decomposition 
$(M,t^0)=\prod_{k=1}^l (M_k,t^{0,k})$ as a product
of germs of F-manifolds.

\medskip
(b) If $E$ is an Euler field of $M$, then $E$
decomposes as $E=\sum_{k=1}^l E_k$ with $E_k$
(the canonical lift of) an Euler field on $M_k$. 
\end{theorem}

\begin{remark}\label{t2.4}\cite[Proposition 2.10]{He02}
Vice versa, if one has $l$ F-manifolds
$(M_k,\circ_k,e_k)$, their product $M=\prod_{k=1}^l M_k$
inherits a natural structure as a germ of an F-manifold.
And if there are Euler fields $E_k$, then the sum 
$E=\sum_{k=1}^l(\textup{lift of }E_k\textup{ to }M)$
is an Euler field on the product $M$.
\end{remark}

Because of Theorem \ref{t2.3}, the first task in a
classification of F-manifolds is the classification 
of irreducible 
germs of F-manifolds, without and with Euler fields.

\begin{lemma}\label{t2.5}\cite[Example 2.12 (i)]{He02}
In dimension 1, (up to isomorphism) there is only
one germ of an F-manifold, the germ $(M,0)=(\C,0)$ with
$e=\paa/\paa_{u_1}$, where $u_1$ is the coordinate on $\C$. 
Any Euler field on it has the shape 
$E=(u_1+c_1)e$ for some $c_1\in\C$. 
\end{lemma}

\begin{definition}\label{t2.6}
(a) A germ of an F-manifold is {\it semisimple}
if it is isomorphic to a product of 1-dimensional germs
of F-manifolds.

(b) An F-manifold is {\it generically semisimple} if it is
semisimple at generic points.
\end{definition}

\begin{remarks}\label{t2.7}
(i) By Lemma \ref{t2.5}, a semisimple germ is isomorphic to 
$(\C^n,0)$ with coordinates $u=(u_1,...,u_n)$ and partial units
$e_k=\paa_{u_k}$, which determine the multiplication by
$e_k\circ e_k=e_k$ and $e_{k_1}\circ e_{k_2}=0$
for $k_1\neq k_2$. The global unit field is $e=\sum_{k=1}^ne_k$.
Any Euler field on this F-manifold has the shape 
$E=\sum_{k=1}^n (u_k+c_k)e_k$ for some $c_1,...,c_n\in\C$.
The semisimple germ of dimension $n$ is said to be of type 
$A_1^n$.

\medskip
(ii) In (i), the coordinates $u_k$ are up to a shift
equal to the eigenvalues $u_k+c_k$ of $E\circ$.
One can also use these eigenvalues as new coordinates.
Then one obtains a germ isomorphic to $(\C^n,-(c_1,...,c_n))$
at some point $-(c_1,...,c_n)\in\C^n$. These coordinates
are {\it Dubrovin's canonical coordinates}.

\medskip
(iii) In the case of a generically semisimple F-manifold,
the {\it caustic} is the set 
\begin{eqnarray}\label{2.2}
\KK:=\{t\in M\,|\, \textup{the germ }(M,t)\textup{ is not
semisimple}\}
\end{eqnarray}
of points where the multiplication is not semisimple.
It is either empty or an analytic hypersurface
\cite[Proposition 2.6]{He02}. 
The {\it Maxwell stratum} is the set
\begin{eqnarray}\label{2.3}
\KK_2&:=&(\textup{the closure of the set }\\
&&\{t\in M-\KK\,|\, \textup{some eigenvalues of }E\circ
\textup{ coincide}\}).\nonumber
\end{eqnarray}
It is either empty or an analytic hypersurface. 
The union of caustic and Maxwell stratum
is the {\it bifurcation set}
$\KK^{bif}:=\KK\cup\KK_2$. It is the set of points where
the eigenvalues of $E\circ$ are not all different.
These three notions come from singularity theory.
\end{remarks}

The 2-dimensional germs of F-manifolds are classified.

\begin{theorem}\label{t2.8}\cite[Theorem 4.7]{He02}
In dimension 2, (up to isomorphism) the germs of
F-manifolds fall into three types:

\medskip
(a) The semisimple germ (of type $A_1^2$). 
See Remark \ref{t2.7} (i) for it 
and for the Euler fields on it.

\medskip
(b) Irreducible germs, which (i.e. some holomorphic 
representatives of them) are at generic points semisimple.
They form a series $I_2(m)$, $m\in\Z_{\geq 3}$.
The germ of type $I_2(m)$ can be given as follows.
\begin{eqnarray}
(M,0)&=&(\C^2,0)\ \textup{with coordinates }t=(t_1,t_2)
\textup{ and }\paa_k:=\frac{\paa}{\paa t_k},\nonumber\\
e&=&\paa_1,\quad \paa_2\circ\paa_2 =t_2^{m-2}e.\label{2.4}
\end{eqnarray}
Any Euler field takes the shape
\begin{eqnarray}\label{2.5}
E&=& (t_1+c_1)\paa_1 + \frac{2}{m}t_2\paa_2
\quad\textup{for some }c_1\in\C.
\end{eqnarray}

\medskip
(c) An irreducible germ, such that the multiplication is
everywhere irreducible. It is called $\NN_2$, and it
can be given as follows.
\begin{eqnarray}
(M,0)&=&(\C^2,0)\ \textup{with coordinates }t=(t_1,t_2)
\textup{ and }
\paa_k:=\frac{\paa}{\paa t_k},\nonumber\\
e&=&\paa_1,\quad \paa_2\circ\paa_2 =0.\label{2.6}
\end{eqnarray}
Any Euler field takes the shape
\begin{eqnarray}\label{2.7}
E&=& (t_1+c_1)\paa_1 + g(t_2)\paa_2
\quad\textup{for some }c_1\in\C\\
&&\hspace*{3cm}\textup{ and some function }
g(t_2)\in\C\{t_2\}.\nonumber
\end{eqnarray}
\end{theorem}

\begin{remarks}\label{t2.9}
(i) Also in the case $m=2$ \eqref{2.4} gives an F-manifold,
and \eqref{2.5} gives an Euler field on it.
But this F-manifold is just the semisimple one.
Its germ at 0 is not irreducible.
And \eqref{2.5} does not give all Euler fields on this 
F-manifold.

\medskip
(ii) The classification of irreducible 3-dimensional germs 
of F-manifolds is incomplete, even in the generically 
semisimple case. First steps are done in 
\cite[5.5]{He02}.

\medskip
(iii) In the case of a generically semisimple germ 
$((M,t^0),\circ,e)$, the automorphism group 
$\Aut((M,t^0),\circ,e)$ is finite \cite[Theorem 4.14]{He02}.
In the case of a germ of type $I_2(m)$ in the normal
form in \eqref{2.4}, it is cyclic of order $m-2$,
and it is generated by the automorphism
$(t_1,t_2)\mapsto (t_1,e^{2\pi i/(m-2)}t_2)$. 
It respects any Euler field.

\medskip
(iv) The automorphism group of the germ $\NN_2$
is in the normal form \eqref{2.6}
\begin{eqnarray}
(\Aut((M,0),\circ,e)&=&\{(t_1,t_2)\mapsto
(t_1,f(t_2))\,|\, f(t_2)\in\C\{t_2\}\nonumber\\
&& \textup{ \ with }f(0)=0,f'(0)\neq 0\}.\label{2.8}
\end{eqnarray}

\medskip
(v) With such an automorphism, one can put the 
Euler field $E$ in \eqref{2.7}  
into a normal form within a family of normal forms,
which has one discrete and two complex parameters, 
so one can erase the dependence
on the functional parameter $g(t_2)\in\C\{t_2\}$ in
\eqref{2.7}. This will be treated in \cite{DH19-2}.
\end{remarks}

All germs of semisimple F-manifolds and of 2-dimensional
F-manifolds have Euler fields.
But already in dimension 3, there are germs of F-manifolds
with no Euler fields 
\cite[Theorem 5.29 (a), Theorem 5.30 (a)(iv)]{He02}.
In the case of irreducible germs of generically semisimple
F-manifolds, the following holds.
This result is part of a general discussion of
existence of Euler fields on F-manifolds
in \cite[3.2]{He02}.

\begin{theorem}\label{t2.10}\cite[Lemma 3.4]{He02}
Let $M$ be a sufficiently small representative of an
irreducible germ $(M,t^0)$ of a generically semisimple
F-manifold. For any $c\in\C$, there is a unique Euler field 
$E_c$ on $M-\KK$ such that for $t\to t^0$ all eigenvalues of
$E\circ$ tend to $c$. We have $E_c=E_0+c\cdot e$. 
The characteristic polynomial of 
$E_c\circ$ extends continuously to $\KK$ and has the
value $(x-c)^{\dim M}$ at $t^0$. 
\end{theorem}

Here observe that the algebra $(T_{t^0},\circ|_{t^0},e|_{t^0})$
is irreducible, so any endomorphism $a\circ|_{t^0}$
for $a\in T_{t^0}M$ has only one eigenvalue.
Therefore if any Euler field exists on $M$, it
arises by holomorphic extension to $\KK$ of some $E_c$.

The holomorphic extendability of the Euler field $E_c$ 
to $\KK$ depends on the geometry of the {\it analytic spectrum}
of the F-manifold, which will be introduced now.

\begin{definition}\label{t2.11}
Let $(M,\circ,e)$ be a complex manifold of dimension $n$ 
with a multiplication $\circ$ on the holomorphic 
tangent bundle and with a unit field $e$.

\medskip
(a) We need some standard data on $T^*M$:
Let $pr_M:T^*M\to M$ denote the projection.
Let $t=(t_1,...,t_n)$ be local coordinates
on $M$, and define $\paa_k:=\paa/\paa t_k$. 
Let $y=(y_1,...,y_n)$ be the fiber coordinates
on $T^*M$ which correspond to $(\paa_1,...,\paa_n)$.
Then the canonical 1-form
$\alpha$ takes the shape $\sum_{i=1}^ny_i{\rm d}t_i$,
and $\omega={\rm d}\alpha$ is the standard symplectic form.
The Hamilton vector field of 
$f\in\OO_{T^*M}$ is 
\begin{eqnarray}\label{2.9}
H_f = \sum_{k=1}^n \Bigl( \frac{\paa f}{\paa t_k}
\cdot \frac{\paa}{\paa y_k} 
-\frac{\paa f}{\paa y_k}\cdot \frac{\paa}{\paa t_k}\Bigr).
\end{eqnarray}
The Posson bracket $\{.,.\}$ on $\OO_{T^*M}$ is defined by
\begin{eqnarray}\label{2.10}
\{f,g\}&:=& H_f(g)=\omega(H_f,H_g)=-H_g(f).
\end{eqnarray}

\medskip
(b) Define an ideal sheaf $\II_M\subset \OO_{T^*M}$
as follows. We choose coordinates $t_k$ and $y_k$ as in 
part (a) and such that $e_1=\paa_1$. 
Write
\begin{eqnarray}\label{2.11}
\paa_i\circ\paa_j =\sum_{k=1}^n a_{ij}^k\paa_k
\textup{ with }a_{ij}^k\in\OO_M.
\end{eqnarray}
Then 
\begin{eqnarray}\label{2.12}
\II_M:=\bigl(y_1-1,y_iy_j-\sum_{k=1}^na_{ij}^ky_k\bigr)
\subset\OO_{T*M}.
\end{eqnarray}
The {\it analytic spectrum} (or {\it spectral cover})
$L_M:=\textup{Specan}_{\OO_M}(TM,\circ)\subset T^*M$
of $(M,\circ,e)$ 
is as a set the set at which the functions in $\II$ vanish.
It is a complex subspace of $T^*M$ with
complex structure given by $\OO_{L_M}=(\OO_M/\II_M)|_{L_M}$.
\end{definition}

The analytic spectrum $L_M$ determines the multiplication
on $M$ and the possible Euler fields.
The projection $L_M\to M$ is finite and flat of degree $\dim M$.
All this is discussed in \cite[2.2 and 3.2]{He02}.
But the following result was missed in \cite{He02}.

\begin{theorem}\label{t2.12}\cite[2.5 Theorem]{HMT09}
A manifold $(M,\circ,e)$ with multiplication $\circ$ on the
tangent bundle and unit field $e$ is an F-manifold
if and only if $\{\II_M,\II_M\}\subset\II_M$.
\end{theorem}

\begin{remarks}\label{t2.13}
(i) In the case of a manifold with a multiplication
and unit field, such that the multiplication is  
generically semisimple, 
the restriction $L_M|_{M-\KK}$ of $L_M$ to $M-\KK$ 
is obviously smooth with $\dim M$ sheets above $M-\KK$.

Theorem 3.2 in \cite{He02} says that then $L_M$ is
reduced everywhere, so also above 
$L_M\cap pr_{T^*M\to M}^{-1}(\KK)$.

\medskip
(ii) In this situation, $\{\II_M,\II_M\}\subset\II_M$
says that $L_M$ is at smooth points a Lagrange submanifold
of $T^*M$. 

\medskip
(iii) But in the case of a manifold with multiplication
and unit field, such that the multiplication is nowhere
semisimple, the analytic spectrum $L_M$ is nowhere
reduced. Then $\II_M$ is quite different from the
reduced ideal $\sqrt{\II_M}$.
Especially, the conditions
\begin{eqnarray}\label{2.13}
\{\II_M,\II_M\}\subset\II_M
\quad\textup{and}\quad  
\{\sqrt{\II_M},\sqrt{\II_M}\}\subset\sqrt{\II_M}
\end{eqnarray}
do not imply one another. 
The second condition in \eqref{2.13} is equivalent
to the condition that $L^{red}_M$ (the reduced space 
underlying $L_M$) is at smooth points a Lagrange submanifold
of $T^*M$. 

\medskip
(iv) The examples 2.5.2 and 2.5.3 in \cite{HMT09}
are examples of F-manifolds (so 
$\{\II_M,\II_M\}\subset\II_M$ holds) with
$\{\sqrt{\II_M},\sqrt{\II_M}\}\not\subset\sqrt{\II_M}$,
so $L_M^{red}$ is at smooth points not a Lagrange
submanifold of $T^*M$. 
We will see in Theorem \ref{t6.7} that these F-manifolds do not
carry $(T)$-structures and therefore cannot be enriched
to flat F-manifolds (or even Frobenius manifolds).

\medskip
(v) The other way round, there are manifolds $(M,\circ,e)$
with multiplication on the tangent bundle and unit field,
which are not F-manifolds
(so $\{\II_M,\II_M\}\not\subset\II_M$ holds),
but which satisfy
$\{\sqrt{\II_M},\sqrt{\II_M}\}\subset\sqrt{\II_M}$,
so $L_M$ is at smooth points a Lagrange submanifold of $T^*M$.
One example with $n=4$ is as follows. $M=\C^4$, 
\begin{eqnarray}\label{2.14}
\II_M&=&(y_1-1,(y_2-t_4y_3)^2,(y_2-t_4y_3)y_3,
y_3^3,y_4-y_3^2),\\
e&=&\paa_1,(\paa_2-t_4\paa_3)^{\circ 2}=0, 
(\paa_2-t_4\paa_3)\circ \paa_3=0,\paa_3^{\circ 3}=0,
\paa_3^{\circ 2}=\paa_4,\nonumber\\
\sqrt{\II_M}&=& (y_1-1,y_2,y_3,y_4),\nonumber 
\end{eqnarray}
\begin{eqnarray*}
H_{(y_2-t_4y_3)y_3}&=& -y_3\frac{\paa}{\paa t_2} 
+(-y_2+2t_4y_3)\frac{\paa}{\paa t_3}-y_3^2\frac{\paa}{\paa y_4},\nonumber\\
\{(y_2-t_4y_3)y_3,y_4-y_3^2\}
&=& H_{(y_2-t_4y_3)y_3}(y_4-y_3^2)=-y_3^2\notin\II_M,\nonumber\\
\textup{so }\{\II_M,\II_M\}&\not\subset&\II_M.\nonumber
\end{eqnarray*}
The second line gives the multiplication in a slightly 
implicit way. For each $t\in M$ the algebra 
$(T_tM,\circ|_t,e|_t)$ is irreducible. Obviously 
$\{\sqrt{\II_M},\sqrt{\II_M}\}\subset\sqrt{\II_M}$ holds. 

\medskip
(vi) We are interested mainly in F-manifolds which
carry $(T)$-structures and can be enriched to
flat F-manifolds (or even Frobenius manifolds).
They must satisfy both conditions in \eqref{2.13}.
This follows from Theorem \ref{t2.12} and Theorem \ref{t6.7}.
But it is not clear how much stronger the 
existence of $(T)$-structures is than both conditions
in \eqref{2.13}. Conjecture \ref{t7.2} says that any
irreducible germ of a generically semisimple
F-manifold with Euler field underlies $(TE)$-structures, 
and it says how many.
\end{remarks}

\begin{remark}\label{t2.14}
Although the classification of F-manifolds is far from
complete, in this paper we take the point of view
that F-manifolds can be generated fairly easily.
One just has to provide in some way an analytic spectrum.
The sections 2.2, 4.5, 5.1, 5.2 and 5.5 in \cite{He02}
give tools for that. In this paper we are mainly 
interested in richer structures above them,
$(TE)$-structures or flat F-manifolds, in obstructions
for their existence and freedom in their choice.
Here we will concentrate mainly, but not only, on
irreducible germs of generically semisimple 
F-manifolds. 
\end{remark}

A generalization of the generically semisimple F-manifolds
are the generically regular F-manifolds.

\begin{definition}\label{t2.15}\cite[Definition 1.2]{DH17}
Let $(M,\circ,e,E)$ be an F-manifold with Euler field.

\medskip
(a) The Euler field is regular at a point $t\in M$
if $E\circ|_t:T_tM\to T_tM$ is a regular endomorphism,
i.e. it has for each eigenvalue only one Jordan block.

\medskip
(b) The F-manifold with Euler field $(M,\circ,e,E)$
is called a [generically] regular F-manifold 
if the Euler field is regular at all
[respectively at generic] points.
\end{definition}

Theorem 1.3 in \cite{DH17} provides a generalization
of the canonical coordinates of a semisimple
F-manifold with Euler field to the case of a 
regular F-manifold.

We expect that some of the results developed in the 
later sections in this paper for generically semisimple
F-manifolds with Euler fields can be generalized to 
generically regular F-manifolds.

\begin{remarks}\label{t2.16}
(i) In order to be able to enrich an F-manifold to a Frobenius manifold,
the algebras $(T_tM,\circ|_t,e|_t)$ must be especially 
{\it Frobenius algebras}.

\medskip
(ii) A finite dimensional (commutative and associative) 
$\C$-algebra with unit is a {\it Frobenius algebra}
if and only if a $\C$-bilinear symmetric non-degenerate
pairing $g$ on $A$ with 
$g(a\circ b,c)=g(a,b\circ c)$ for $a,b,c\in A$ exists.
Such a pairing is called a 
{\it multiplication invariant metric}.

\medskip
(iii) If $A=\bigoplus_{k=1}^l A_k$ is the decomposition
into irreducible algebras (with $A_{k_1}\circ A_{k_2}=0$
for $k_1\neq k_2$), then a multiplication invariant
metric on $A$ restricts to multiplication invariant
metrics on the subalgebras $A_k$. 
And an irreducible algebra $A$ is a Frobenius algebra 
if and only if its socle $\textup{Ann}_A({\bf m}_A)$
has dimension 1. Here ${\bf m}_A$ is the maximal
ideal in $A$.

\medskip
(iv) It is easy to give examples of F-manifolds $M$
whose algebras $(T_tM,\circ|_t,e|_t)$ are Frobenius
algebras only outside some hypersurface $H$.
Then a possible Frobenius manifold structure on $M-H$
cannot extend to $H$.
\end{remarks}

\begin{remarks}\label{t2.17}
These last remarks in this section give a panorama
of F-manifolds.

\medskip
(i) Probably the most important nontrivial F-manifolds
are those F-manifolds $M$ whose analytic spectrum $L_M$
is smooth. They are generically semisimple, but
at points of the caustic, they become interesting.
In \cite[Theorem 5.6]{He02} results of H\"ormander
and Arnold on Lagrange maps are translated into a
1:1-correspondence between irreducible germs $(M,t^0)$
of generically semisimple F-manifolds with smooth
analytic spectrum and stable right equivalence classes
of holomorphic function germs $f:(\C^{n+1},0)\to (\C,0)$
with isolated singularity at 0.
There the germ $(M,t^0)$ is the base space of a
universal unfolding of the hypersurface singularity.
These objects have been studied a lot in singularity theory.

\medskip
(ii) A related family is the family of irreducible germs 
$(M,t^0)$ of generically semisimple F-manifolds whose analytic
spectrum $L_M$ is a union of two smooth components which
intersect transversely in codimension 1.
There is a 1:1-correspondence between such germs and
stable right equivalence classes of {\it boundary 
singularities} \cite[Theorem 5.14]{He02}.
This builds on work of Nguyen tien Dai and Nguyen huu Duc.
Again, $(M,t^0)$ arises as the base space of a universal
unfolding of a boundary singularity.

\medskip
(iii) Long before the more general work in \cite{KMS18},
K. Saito \cite{SaK79} established on the complexified
orbit spaces of the finite Coxeter groups the structure
of a Frobenius manifold with Euler field, where
all data are polynomial. This includes an F-manifold
structure. Dubrovin \cite{Du96} recovered these
Frobenius manifolds. The second author studied the
underlying F-manifolds and related F-manifolds 
and proved a conjecture of Dubrovin: 
A germ of a generically semisimple
Frobenius manifold with an Euler field with positive weights
decomposes uniquely into a product of germs at 0 of 
Frobenius manifolds from finite Coxeter groups
\cite[Theorem 5.25]{He02}

\medskip
(iv) Subsection 5.5 in \cite{He02}
gives several families of germs $(M,0)$ of generically
semisimple 3-dimensional F-manifolds, with and without 
Euler fields, some with $(T_0M,\circ|_0,e|_0)
\cong \C[x]/(x^3)$ (a Frobenius algebra),
some with $(T_0M,\circ|_0,e|_0)
\cong \C[x,y]/(x^2,xy,y^2)$ (not a Frobenius algebra).
These two algebras are the only irreducible 3-dimensional 
(commutative) $\C$-algebras.
Theorem 5.30 gives a complete classification of
irreducible germs of 3-dimensional F-manifolds
with $(T_0M,\circ|_0,e|_0)\cong \C[x]/(x^3)$ and 
whose analytic spectrum has 3 components. 

\medskip
(v) Givental's paper on Lagrange maps \cite{Gi88}
contains implicitly (via their analytic spectra)
many families of generically semisimple F-manifolds
\cite[Remark 5.33]{He02}.
\end{remarks}

\section{Frobenius manifolds and slightly weaker structures}
\label{c3}
\setcounter{equation}{0}

\noindent 
Dubrovin defined Frobenius manifolds in \cite{Du92}
and studied them in \cite{Du96}, \cite{Du99}, \cite{Du04} and
many other papers. They arise in integrable systems,
in quantum cohomology and in singularity theory.
One version of mirror symmetry compares Frobenius
manifolds from quantum cohomology and from singularity
theory. 

In this section we define them together with slightly
weaker structures, which had been defined
by Manin \cite{Ma05} long ago, but which received
some new attention because of relations to Painlev\'e
equations and finite complex reflection groups.

We will say something about recent work of three groups, 
Arsie \& Lorenzoni, Kato \& Kawakami \& Mano \& Sekiguchi, 
Konishi \& Minabe \& Shiraishi.

\begin{definition}\label{t3.1}
(a) \cite{Ma05}
A {\it flat F-manifold} $(M,\circ,e,D)$ is an
F-manifold $(M,\circ,e)$ together with a flat and torsion free 
connection $D$ on the holomorphic tangent bundle $TM$ 
which satisfies
\begin{eqnarray}\label{3.1}
D(C^M)&=& 0,\quad D(e)=0.
\end{eqnarray}
Here $C^M$ is the Higgs field on $TM$ from the 
multiplication, i.e
\begin{eqnarray}\label{3.2}
C^M_X=X\circ :TM\to TM \quad\textup{for }X\in\TT_M.
\end{eqnarray}
More explicitly, $D(C^M)=0$ says 
\begin{eqnarray}\label{3.3}
D_X(Y\circ)-D_Y(X\circ) &=& [X,Y]\circ\quad\textup{for }
X,Y\in\TT_M.
\end{eqnarray}

(b) A {\it flat F-manifold with Euler field} $(M,\circ,e,E,D)$
is a flat F-manifold $(M,\circ,e,D)$ together with an
Euler field $E$ of the F-manifold such that 
$D_\bullet E:TM\to TM,$ $X\mapsto D_XE$, 
is a flat endomorphism, i.e.
$D(D_\bullet E)=0$.

\medskip
(c) \cite{Du92}\cite{Du96}
A {\it Frobenius manifold} $(M,\circ,e,g)$
is an F-manifold $(M,\circ,e)$ together with a
holomorphic bilinear symmetric nondegenerate pairing $g$ 
(called {\it metric}) on $TM$ which is 
{\it multiplication invariant}, i.e. 
\begin{eqnarray}\label{3.4}
g(X\circ Y,Z)=g(X,Y\circ Z)\quad\textup{for }X,Y,Z\in\TT_M,
\end{eqnarray}
and such that its Levi-Civita connection $D$ together
with the F-manifold is a flat F-manifold $(M,\circ,e,D)$.

\medskip
(d) A {\it Frobenius manifold with Euler field}
$(M,\circ,e,E,g)$ is a Frobenius manifold
$(M,\circ,e,g)$ together with an Euler field
$E$ of the F-manifold such that 
\begin{eqnarray}\label{3.5}
\Lie_E(g)=(2-d)\cdot g\quad\textup{for some }d\in\C.
\end{eqnarray}
\end{definition}

\begin{remarks}\label{t3.2}
(i) If $(M,\circ,e)$ is any manifold with multiplication
on the tangent bundle and unit field, and if $D$ is 
any torsion free connection on $TM$ with $D(C^M)=0$,
then $(M,\circ,e)$ is an F-manifold
\cite[Theorem 2.14]{He02}. Therefore the definitions
above are slightly redundant. 
The condition $D(C^M)=0$ is called {\it potentiality}
because of (iv) and (v) below.

\medskip
(ii) In the presence of a metric $g$ with \eqref{3.4},
one has even an equivalence \cite[Theorem 2.15]{He02}:
Let $(M,\circ,e,g)$ be a manifold with multiplication
on the tangent bundle and unit field and metric $g$
(so $g$ is a holomorphic bilinear symmetric nondegenerate
pairing) with \eqref{3.4}, 
and let $D$ be the Levi-Civita connection of $g$.
Then $D(C^M)=0$ is equivalent to $(M,\circ,e)$
being an F-manifold and the 1-form $g(.,e)$ being closed.

\medskip
(iii) These definitions do not stress the role of the
$D$-flat vector fields, therefore those in (c) and (d)
look different from Dubrovin's definitions.
But they are equivalent. 

\medskip
(iv) Locally in a flat F-manifold one can choose 
coordinates $t=(t_1,...,t_n)$
with $D$-flat coordinate vector fields $\paa_k=\paa/\paa t_k$
and with $e=\paa_1$. Let $a_{ij}^k\in\OO_M$ be the 
coefficients of the multiplication with respect to these
vector fields as in \eqref{2.11}. Then the
potentiality $D(C^M)=0$ gives
\begin{eqnarray}\label{3.6}
\paa_i(a_{jk}^l) = \paa_j(a_{ik}^l),
\end{eqnarray}
so locally there are functions $b_k^l\in \OO_M$ with
$\paa_i(b_k^l)=a_{ik}^l$. Now the symmetry 
\begin{eqnarray}\label{3.7}
\paa_i(b_k^l) = a_{ik}^l = a_{ki}^l =\paa_k(b_i^l)
\end{eqnarray}
gives locally a function $c^l\in\OO_M$ with
$\paa_ic^l=b_i^l$. The vector field
$\sum_lc^l\paa_l$ is called a {\it vector potential}
in \cite{Ma05}. It satisfies 
$\paa_i\circ\paa_j=[\paa_i,[\paa_j,\sum_lc^l\paa_l]]$. 
This shows the equivalence of the
definition above of a flat F-manifold with 
Definition 2.2 b) in \cite{Ma05}.

\medskip
(v) In the case of a Frobenius manifold one can integrate
once more. $\paa_i\paa_jc^l=a_{ij}^l$ and \eqref{3.4} give
locally a {\it potential} $F\in\OO_M$ with
\begin{eqnarray}\label{3.8}
\paa_i\paa_j\paa_k F &=& g(\paa_i\circ\paa_j,\paa_k).
\end{eqnarray}
In the case of quantum cohomology, this potential
is the Gromov-Witten potential of genus 0 Gromov-Witten 
invariants (without gravitational descendents), 
and there it is the starting point for the construction 
of the Frobenius manifold.

\medskip
(vi) But in singularity theory, one first has the
F-manifold, and then one constructs a flat metric
and its Levi-Civita connection. 
A general frame for this construction will be recalled 
below in Theorem \ref{t6.6}. 
The potential $F$ from (v) or the vector potential 
$\sum_l c^l\paa_l$ from (iv) do not play any role 
in this approach to Frobenius manifolds and flat F-manifolds. 
\end{remarks}

\begin{remarks}\label{t3.3}
(i) Dubrovin's paper \cite{Du04} starts with the observation
that a Frobenius manifold with Euler field $(M,\circ,e,E,g)$
gives rise to a second structure 
$(M-\DD,\ast,E,e,\eta)$ which is almost a Frobenius manifold: Here $\DD=\{t\in M\,|\, \det E\circ =0\}$ is the 
discriminant (either empty or a hypersurface), 
$\ast$ is a new multiplication $\ast$ with 
$X\ast Y:= E^{-1}\circ X\circ Y$ and unit field $E$, 
and $\eta$ is a new metric with 
$\eta(X,Y)=g(E^{-1}\circ X, Y)$. Then
$(M-\DD,\ast,E)$ satisfies all properties of a 
Frobenius manifold, except that the new unit field $E$ 
and the Levi-Civita connection $D^{(\eta)}$ of $\eta$
satisfy $D^{(\eta)}_\bullet E=\frac{1-d}{2}\id$,
so for $d\neq 1$, $E$ is not $D^{(\eta)}$-flat.
The old unit field $e$ satisfies compatibilities with
$\ast$ and $\eta$ different from those of an Euler field.

\medskip
(ii) Manin \cite{Ma05} extended this duality to the case
of flat F-manifolds. On the dual side he had 1 parameter
and a 1-parameter family of connections.

\medskip
(iii) David and Strachan \cite{DS11} extended
this duality further by replacing the Euler field by an
{\it eventual identity} and studied several enrichments
\cite{DS14}.

\medskip
(iv) Arsie and Lorenzoni \cite{Lo14}
considered both structures together, observed additional
joint compatibility conditions, and called the total
structure {\it bi-flat F-manifold}.

\medskip
(v) Konishi, Minabe and Shiraishi \cite{KMS18}
put this into a good conceptual framework.
Their Proposition 3.1 and Definition 3.5 allow to identify
the 1-parameter family of connections on the dual side
as the 1-parameter family of second structure connections
in \cite[(9.10) in Definition 9.3]{He02} 
with $z=0$ and with parameter $s$.

Also, they showed in Lemma 4.2 and Lemma 4.3 in \cite{KMS18}
that and how a bi-flat F-manifold is equivalent to a
flat F-manifold with Euler field. 

\medskip
(v) Sabbah calls in \cite[VII Definition 1.1]{Sa02}
a flat F-manifold with Euler field a 
{\it Saito structure without metric}.
\cite{KMS18} and [KM19] follow this notation.
In \cite[VII Proposition 2.2]{Sa02} he states that
a Saito structure with metric is indeed a Frobenius
manifold with Euler field.
\end{remarks}

\begin{remarks}\label{t3.4}
(i) Kato, Mano and Sekiguchi \cite{KMS15}, 
Arsie and Lorenzoni \cite{AL17}
and Konishi, Minabe and Shiraishi \cite{KMS18}
all establish and study the structures of flat 
F-manifolds with Euler field on the orbit spaces of most 
finite complex reflection groups.

\cite{KMS15} and \cite{KMS18} establish them on all
finite well-generated complex reflection groups,
\cite{AL17} restricts to rank 2 and rank 3 groups.
All three papers are explicit about many rank 2 groups.

\medskip
(ii)  Before their papers, in the case of the F-manifolds
$I_2(m)$ ($m\geq 3$) a 1-parameter family of flat
F-manifolds with Euler fields  was known, 
which contains as one member a Frobenius manifold
with Euler field.
\cite{KMS15} and \cite{AL17} found independently in the
cases $I_2(m)$ with $m$ even that this is part of 
2-parameter family, and that the flat F-manifolds
with Euler fields which are canonically associated
to most of the finite rank 2 complex reflection groups
are members of the 2-parameter family, which are not in
the old 1-parameter family. 

\medskip
(iii) Theorem \ref{t6.6}, Theorem \ref{t8.5} and Remark
\ref{t8.8} below will show that the old 1-parameter family
for odd $m$ and the new 2-parameter family for even $m$
comprise all flat F-manifolds with Euler fields over
the F-manifold $I_2(m)$.
Conjecture \ref{t7.2} puts these families into a 
general frame.
\end{remarks}

\begin{remarks}\label{t3.5}
(i) Dubrovin \cite{Du96} had observed a relation
between 3-dimensional semisimple Frobenius manifolds
with Euler field and the solutions of a 1-parameter
subfamily of the (4-dimensional family of the)
Painlev\'e VI equations.

\medskip
(ii) Kato, Mano and Sekiguchi \cite[Corollary 4.16]{KMS15} 
and Lorenzoni \cite[Theorem 4.1]{Lo14} found independently 
a generalization.
They observed a relation between 3-dimensional semisimple 
flat F-manifolds with Euler field and 
generic solutions of all Painlev\'e VI equations. 

\medskip
(iii) Arsie and Lorenzoni 
\cite[Theorem 8.7 and Theorem 8.4]{AL15}
and Kawakami and Mano \cite[Theorem 4.6]{KM19}
generalized this further to relations 
between 3-dimensional regular flat F-manifolds
with Euler fields and generic solutions of all
Painlev\'e equations of types VI, V and IV.
And Kawakami and Mano \cite[Theorem 5.4]{KM19}
also got relations between
4-dimensional regular flat F-manifolds with
Euler fields and generic solutions of all
Painlev\'e equations of types VI, V, IV, III and II.
The following table gives more information.
Here a tuple $(a_1,...,a_l)\in\N^l$ means
that regular F-manifolds are considered where
$E\circ$ has at each point $l$ Jordan blocks
of sizes $a_1,...,a_l$ with different eigenvalues.

\medskip
\begin{tabular}{cccc}
 & regular flat F-manifolds & & solutions of Painlev\'e\\ 
 & with Euler fields & $\leftrightarrow$ & equations\\ \hline
dimension &  size of Jordan blocks & & type of the equations \\
\hline 
3 & (1,1,1) & $\leftrightarrow$ & VI \\
3 & (2,1) & $\leftrightarrow$ & V \\
3 & (3) & $\leftrightarrow$ & IV \\ \hline
4 & (1,1,1,1) & $\leftrightarrow$ & VI \\
4 & (2,1,1) & $\leftrightarrow$ & V \\
4 & (3,1) & $\leftrightarrow$ & IV \\
4 & (2,2) & $\leftrightarrow$ & III \\
4 & (4) & $\leftrightarrow$ & II
\end{tabular}

\medskip
Though it is not claimed that the correspondences
are 1:1. \cite{KM19} speak of correspondences
between generic objects on both sides.
\cite{Lo14} and \cite{AL15} say that the objects
on the left hand side are parametrized locally by 
solutions of the full Painlev\'e equations.

\medskip
(iv) Arsie and Lorenzoni consider on the left hand
side bi-flat F-manifolds. Due to \cite[Lemma 4.2 and Lemma
4.3]{KMS18} these are equivalent to flat F-manifolds
with Euler fields (Remark \ref{t3.3} (v)).

\medskip
(v) The correspondences in the table above
are beautiful. They motivate a study of 
3-dimensional (and 4-dimensional) F-manifolds beyond
\cite[5.5]{He02}.
\end{remarks}

\section{A dictionary of connections with different
enrichments}\label{c4}
\setcounter{equation}{0}

Definition \ref{t3.1} considers 4 cases,
flat F-manifolds without/with Euler field and
without/with metric. 
Theorem \ref{t6.6} will recall a well-known recipe
for the construction of these 4 structures.

This section defines and discusses structures which
arise in this recipe, holomorphic vector bundles
with meromorphic connections with certain properties
and certain enrichments.
This section can be seen as a dictionary of these structures.
It gives only definitions and elementary consequences.
The structures had been considered before in \cite{HM04},
and they are related to structures in \cite[VII]{Sa02}
and in \cite{Sa05}.

\begin{definition}\label{t4.1}
(a) Definition of a {\it (T)-structure} 
$(H\to\C\times M,\nabla)$:
$H\to\C\times M$ is a holomorphic vector bundle. 
$\nabla$ is a map
\begin{eqnarray}\label{4.1}
\nabla: \OO(H)\to \frac{1}{z}\OO_{\C\times M}\cdot \Omega^1_M\otimes \OO(H),
\end{eqnarray}
which satisfies the Leibniz rule,
$$\nnn_X(a\cdot s)= X(a)\cdot s+a\cdot \nnn_Xa
\quad\textup{for }X\in\TT_M,a\in\OO_{\C\times M},s\in \OO(H),$$
and which is flat with respect to $X\in\TT_M$,
$$\nnn_X\nnn_Y-\nnn_Y\nnn_X=\nnn_{[X,Y]}
\quad\textup{for }X,Y\in \TT_M.$$
Equivalent: For any $z\in\C^*$, the restriction of $\nabla$ to 
$H|_{\{z\}\times M}$ is a flat holomorphic connection.

\medskip
(b) Definition of a {\it (TE)-structure} 
$(H\to\C\times M,\nabla)$:
$H\to\C\times M$ is a holomorphic vector bundle. 
$\nabla$ is a flat connection on $H|_{\C^*\times M}$ with a pole
of Poincar\'e rank 1 along $\{0\}\times M$, so it is a map
\begin{eqnarray}\label{4.2}
\nabla: \OO(H)\to 
\bigl(\frac{1}{z}\OO_{\C\times M}\cdot\Omega^1_M
+\frac{1}{z^2}\OO_{\C\times M}\cdot{\rm d}z\bigr)\otimes\OO(H)
\end{eqnarray}
which satisfies the Leibniz rule and is flat.

\medskip
(c) Definition of a {\it (TP)-structure 
$(H\to\C\times M,\nabla,m,P)$ of weight $m\in\Z$}:
Denote by $j$ the holomorphic involution
\begin{eqnarray}\label{4.3}
j:\P^1\times M\to\P^1\times M,\ (z,t)\mapsto (-z,t).
\end{eqnarray}
A ($TP)$-structure of weight $m$ 
is a $(T)$-structure together with an
$\OO_{\C\times M}$-bilinear $(-1)^m$-symmetric $\nnn$-flat
pairing 
\begin{eqnarray}\label{4.4}
P: \OO(H) \times j^*\OO(H) \to z^{m}\OO_{\C\times M},
\end{eqnarray}
which is nondegenerate in the following sense:
$z^{-m}\cdot P$ is nondegenerate at each $(z,t)\in\C\times M$.
Pointwise $P$ is given by a $\C$-bilinear form
\begin{eqnarray}\label{4.5}
P: H_{z,t}\times H_{-z,t}\to\C \quad\textup{for }
(z,t)\in \C^*\times M.
\end{eqnarray}
which is nondegenerate and $(-1)^m$-symmetric.

\medskip
(d) Definition of a {\it (TEP)-structure} 
$(H\to\C\times M,\nabla,m,P)$ of weight $m\in\Z$:
It is simultaneously a $(TE)$-structure and a $(TP)$-structure
where $P$ is $\nabla$-flat.

\medskip
(e) Definition of a {\it (TL)-structure} $(H\to\P^1\times M,\nabla)$:
$H\to\P^1\times M$ is a holomorphic vector bundle. 
$\nabla$ is a map
\begin{eqnarray}\label{4.6}
\nabla: \OO(H)\to \left(\frac{1}{z}\OO_{\P^1\times M}+\OO_{\P^1\times M}\right)
\cdot \Omega^1_M\otimes \OO(H),
\end{eqnarray}
such that for any $z\in\P^1-\{0\}$, the restriction 
of $\nabla$ to  $H|_{\{z\}\times M}$ is a flat connection.
It is called {\it pure} if for any $t\in M$ the restriction
$H|_{\P^1\times\{t\}}$ is a trivial holomorphic bundle on $\P^1$.

\medskip
(f) Definition of a {\it (TLE)-structure} 
$(H\to\P^1\times M,\nabla)$:
It is simultaneously a $(TE)$-structure and a 
$(TL)$-structure, where the connection $\nabla$ 
has a logarithmic pole along $\{\infty\}\times M$.
The $(TLE)$-structure is called {\it pure} if the
$(TL)$-structure is pure.

\medskip
(g)  Definition of a {\it (TLP)-structure} 
$(H\to\P^1\times M,\nabla,m,P)$ of weight $m\in\Z$:
It is simultaneously a $(TL)$-structure and a 
$(TP)$-structure where the pairing $z^{-m}\cdot P$ extends to a holomorphic symmetric 
and everywhere nondegenerate pairing
\begin{eqnarray}\label{4.7}
z^{-m}\cdot P: \OO(H)\times j^*\OO(H)\to \OO_{\P^1\times M}.
\end{eqnarray}
The $(TLP)$-structure is called {\it pure} if the
$(TL)$-structure is pure. 

\medskip
(h) Definition of a {\it (TLEP)-structure} 
$(H\to\P^1\times M,\nabla,m,P)$ of weight $m\in\Z$:
It is simultaneously a $(TEP)$-structure,
a $(TLE)$-structure and a $(TLP)$-structure.
It is called {\it pure} if the $(TL)$-structure is pure. 
\end{definition}

\begin{remark}\label{t4.2}
Here we fix some useful conventions how to represent
endomorphisms and bilinear forms by matrices.
An endomorphism $f:V\to V$ of a vector space
$V$ of dimension $r\in\N$ 
is after the choice of a basis $\uuuu{v}=(v_1,...,v_r)$
of the vector space represented by the matrix 
$f^{mat}=(f^{mat}_{ij})\in M_{r\times r}(\C)$ with
\begin{eqnarray}\label{4.8}
f(\uuuu{v})=\uuuu{v}\cdot f^{mat}:=(\sum_{i=1}^rf^{mat}_{i1}
\cdot  v_i,... ,\sum_{i=1}^r f^{mat}_{ir}\cdot v_i).
\end{eqnarray}
This way (contrary to the transpose way) has the advantage
that for two endomorphisms $f$ and $g$ 
$(f\circ g)^{mat}=f^{mat}\cdot g^{mat}$.
And an abstract vector $\sum_{i=1}^rx_iv_i=\uuuu{v}\cdot x$
with $x=(x_1...x_r)^t$ 
is represented in $M_{r\times 1}(\C)$ by this column vector $x$,
and $f^{mat}\cdot x$ represents $f(\uuuu{v}\cdot x)$.

A bilinear form $h:V\times V$ is represented by the matrix
\begin{eqnarray}\label{4.9}
h(\uuuu{v}^t,\uuuu{v}):=(h(v_i,v_j)).
\end{eqnarray}
Then $h(f(\uuuu{v})^t,g(\uuuu{v}))=(f^{mat})^t\cdot 
h(\uuuu{v}^t,\uuuu{v})\cdot g^{mat}$.
These conventions are for example also used in \cite{Ko87}.
\end{remark}

\begin{remarks}\label{t4.3}
Here we write the data in Definition \ref{t4.1} (a)--(d) 
and the compatibility conditions between them in terms
of matrices. 
Consider a $(TEP)$-structure $(H\to\C\times M,\nnn,m,P)$
of weight $m\in\Z$ and rank $\rk H=r\in\N$. 
We will fix the notations for a trivialization 
of the bundle $H|_{V\times M}$ for some small neighborhood
$V\subset \C$ of 0. Trivialization means the choice of a
basis $\uuuu{v}=(v_1,...,v_r)$ of the bundle $H|_{V\times M}$.
Choose local coordinates $t=(t_1,...,t_n)$ 
with coordinate vector fields $\paa_i=\paa/\paa t_i$ on $M$.
We write
\begin{eqnarray}
\nnn\uuuu{v}&=&\uuuu{v}\cdot\Omega\quad\textup{with}\nonumber\\
\Omega &=& \sum_{i=1}^r z^{-1}\cdot A_i(z,t)\ddd t_i
+z^{-2}B(z,t)\ddd z,\label{4.10}\\
A_i(z,t)&=& \sum_{k\geq 0}A_i^{(k)}z^k\in
M_{r\times r}(\OO_{V\times M}),
\label{4.11}\\
B(z,t)&=&\sum_{k\geq 0} B^{(k)}z^k\in
M_{r\times r}(\OO_{V\times M}),\label{4.12}
\end{eqnarray}
and 
\begin{eqnarray}
z^{-m}\cdot P(\uuuu{v}^t,\uuuu{v})&=& P^{mat}(z,t)
= \sum_{k\geq 0}P^{mat,(k)}z^k\in M_{r\times r}
(\OO_{V\times M})\hspace*{1cm}\label{4.13}\\
\textup{with }(P^{mat,(k)})^t&=& (-1)^k\cdot P^{mat,(k)},\label{4.14}
\end{eqnarray}
and $A_i^{(k)},B^{(k)},P^{mat,(k)}\in\OO_M$, 
but this dependence on $t\in M$ is usually not written
explicity. The matrix $P^{mat,(0)}$ is by hypothesis 
nondegenerate at each $t\in M$. 

The flatness $0=\ddd \Omega+\Omega\land\Omega$ of the connection
$\nnn$ says for $i,j\in\{1,...,n\}$ with $i\neq j$
\begin{eqnarray}\label{4.15}
0&=& z\paa_iA_j-z\paa_jA_i+[A_i,A_j],\\
0&=& z\paa_i B-z^2\paa_z A_i + zA_i + [A_i,B].\label{4.16}
\end{eqnarray}
These equations split into the parts for the different powers
$z^k$ for $k\geq 0$ as follows (with $A_i^{(-1)}=B^{(-1)}=0$),
\begin{eqnarray}\label{4.17}
0&=& \paa_iA_j^{(k-1)}-\paa_jA_i^{(k-1)}
+\sum_{l=0}^k[A_i^{(l)},A_j^{(k-l)}],\\
0&=& \paa_i B^{(k-1)}-(k-2)A_i^{(k-1)} 
+\sum_{l=0}^k[A_i^{(l)},B^{(k-l)}].\label{4.18}
\end{eqnarray}
The flatness of $P$ is expressed by 
\begin{eqnarray}
z\paa_i P^{mat}(z,t)&=& A_i^t(z,t)P^{mat}(z,t)
-P^{mat}(z,t)A_i(-z,t),\label{4.19}\\
z^2\paa_z P^{mat}(z,t)&=& -m\cdot zP^{mat}(z,t)\label{4.20}\\
&& +B^t(z,t) P^{mat}(z,t)-P^{mat}(z,t)B(-z,t).\nonumber
\end{eqnarray}
These equations split into the parts for the different 
powers $z^k$ for $k\geq 0$ as follows (with $P^{mat,(-1)}=0$),
\begin{eqnarray}\label{4.21}
&&\paa_i P^{mat,(k-1)}=\sum_{l=0}^k \Bigl((A_i^{(l)})^tP^{mat,(k-l)}
-(-1)^l P^{mat(k-l)}A_i^{(l)}\Bigr),\hspace*{1cm}\\
&&(k-1+m)P^{mat,(k-1)}= \sum_{l=0}^k \Bigl((B^{(l)})^tP^{mat,(k-l)}\Bigr.
\nonumber\\
&&\hspace*{5cm} \Bigl. -(-1)^lP^{mat,(k-l)}B^{(l)}\Bigr).\label{4.22}
\end{eqnarray}

In the case of a $(TE)$-structure, $P^{mat}$ and all
equations which contain it are dropped.
In the case of a $(TP)$-structure, the summand 
$z^{-2}B(t,z)\ddd z$ in $\Omega$ and all equations
except \eqref{4.10} which contain $B$ are dropped.
In the case of a $(T)$-structure, $B$ and $P$ and all
equations except \eqref{4.10} which contain them are dropped. 

Consider a second $(TEP)$-structure 
($\www H\to\C\times M,\www\nnn,m,\www P)$ 
of weight $m\in\Z$ over $M$, 
where all data (except $M$ and $m$) are written with a tilde. 
Let $\uuuu{v}$ and $\uuuu{\www v}$ be trivializations.
A holomorphic isomorphism from the first to the second
$(TEP)$-structure maps $\uuuu{v}\cdot T$ to $\uuuu{\www v}$,
where 
$T=T(z,t) =\sum_{k\geq 0}T^{(k)}z^k\in M_{r\times r}(\OO_{(\C,0)\times M})$
with $T^{(k)}\in M_{r\times r}(\OO_{M})$ and $T^{(0)}$ 
invertible satisfies 
\begin{eqnarray}\label{4.23}
\uuuu{v}\cdot\Omega\cdot T+\uuuu{v}\cdot\ddd T=
\nnn(\uuuu{v}\cdot T)=\uuuu{v}\cdot T\cdot\www\Omega,\\
T^t \cdot P^{mat}\cdot T =
z^{-m} P((\uuuu{v}\cdot T)^t,\uuuu{v}\cdot T)
=z^{-m} \www P(\uuuu{\www v}^t,\uuuu{\www v})={\www P}^{mat}.
\label{4.24}
\end{eqnarray}
\eqref{4.23} says more explicitly
\begin{eqnarray}\label{4.25}
0&=& z\paa_i T+A_i\cdot T-T\cdot\www A_i,\\
0&=& z^2\paa_z T+B\cdot T-T\cdot \www B.\label{4.26}
\end{eqnarray}
These equations split into the parts for the different
powers $z^k$ for $k\geq 0$ as follows (with $T^{(-1)}:=0$):
\begin{eqnarray}\label{4.27}
0= \paa_i T^{(k-1)}+\sum_{l=0}^k (A_i^{(l)}\cdot T^{(k-l)}
-T^{(k-l)}\cdot\www A_i^{(l)}),\\
0= (k-1)T^{(k-1)}+\sum_{l=0}^k(B^{(l)}\cdot T^{(k-l)}
-T^{(k-l)}\cdot \www B^{(l)}),\label{4.28}\\
\sum_{i,j\geq 0:\, i+j\leq k}
(T^{(i)})^t\cdot P^{mat,(k-i-j)}\cdot T^{(j)}={\www P}^{mat,(k)}.
\label{4.29}
\end{eqnarray}
\end{remarks}

Definition \ref{t4.4} (a)--(d) fixes structures which are
induced canonically from the structures in Definition 
\ref{t4.1} (a)--(d). Definition \ref{t4.4} (e)--(h)
fixes structures which are induced after some choice
from the structures in Definition \ref{t4.1} (a)--(d),
see Lemma \ref{t4.5}.

\begin{definition}\label{t4.4}
(a) A {\it Higgs bundle} $(K\to M,C)$ is a holomorphic vector bundle $K\to M$
with an $\OO_M$-linear map 
$C:\OO(K)\to \Omega^1_M\otimes \OO(K)$ with $C\land C=0$, i.e. 
the holomorphic endomorphisms $C_X:K\to K$ and $C_Y:K\to K$ commute,
\begin{eqnarray}\label{4.30}
C_XC_Y=C_YC_X \quad\textup{for }X,Y\in\TT_M.
\end{eqnarray}
$C$ is called a {\it Higgs field}.

\medskip
(b) A {\it Higgs bundle with Euler endomorphism}
$(K\to M,C,\UU)$ is a Higgs bundle $(K\to M,C)$
together with a holomorphic endomorphism $\UU:K\to K$ with
$C\UU=\UU C$, i.e. 
\begin{eqnarray}\label{4.31}
C_X\UU =\UU C_X\quad\textup{for }X\in\TT_M.
\end{eqnarray}

\medskip
(c) A {\it Higgs bundle with metric} $(K\to M,C,g)$
is a Higgs bundle $(K\to M,C)$ with a holomorphic
$\C$-bilinear symmetric and nondegenerate pairing
$g:K_t\times K_t\to\C$ for $t\in M$
with 
\begin{eqnarray}\label{4.32}
g(C_X a,b)= g(a,C_X b)\quad\textup{for }X\in \TT_M,\ 
a,b\in \OO(K).
\end{eqnarray}

\medskip
(d) A {\it Higgs bundle with Euler endomorphism and metric}
$(K\to M,C,\UU ,g)$ is simultaneously a Higgs bundle
with Euler endomorphism and a Higgs bundle with metric,
with the additional compatibility condition
\begin{eqnarray}\label{4.33}
g(\UU a,b)=g(a,\UU b)\quad\textup{for }a,b\in\OO(K).
\end{eqnarray}

\medskip
(e) A {\it Higgs bundle with good connection} $(K\to M,C,D)$ 
is a Higgs bundle $(K\to M,C)$ with a holomorphic connection 
$D$ on $K$ which satisfies the {\it potentiality condition} 
$D(C)=0$, i.e. 
\begin{eqnarray}\label{4.34}
D_X(C_Y)-D_Y(C_X)
=C_{[X,Y]}\quad \textup{for } X,Y\in\TT_M.
\end{eqnarray}

\medskip
(f) A {\it good pair $(\UU ,Q)$ of endomorphisms} on a Higgs 
bundle with good connection $(K\to M,C,D)$ consists of two 
holomorphic endomorphisms $\UU :K\to K$ and $Q:K\to K$ with
$C\UU =\UU C$ (i.e. \eqref{4.31}) and 
\begin{eqnarray}
0= D(\UU )-[C,Q]+C .\label{4.35}
\end{eqnarray}

\medskip
(g) A {\it Higgs bundle with good connection and metric 
and second order pairing}
$(K\to M,C,D,g,g^{(1)})$ is a Higgs bundle with 
good connection $D$ and metric $g$ 
and a holomorphic $\C$-bilinear skew-symmetric pairing 
$g^{(1)}:K_t\times K_t\to\C$ for $t\in M$ with 
\begin{eqnarray}\label{4.36}
D_X(g)(a,b)\bigl(= D_X(g(a,b))-g(D_Xa,b)-g(a,D_Xb)\bigr)\\
=g^{(1)}(C_Xa,b)-g^{(1)}(a,C_Xb)\quad \textup{for }
X\in\TT_M,a,b\in\OO(K).\nonumber
\end{eqnarray}

\medskip
(h) A {\it Higgs bundle with good connection and
good pair of endomorphisms and metric and second order
pairing and weight $m\in\Z$} 
$(K\to M,C,D,\UU ,Q,g,g^{(1)},m)$ is simultaneously
a Higgs bundle with good connection and good pair of
endomorphisms and a Higgs bundle with good connection
and metric and second order pairing, with the following
additional compatibility conditions,
\begin{eqnarray}\label{4.37}
&&g(Q a,b)+g(a,Q b)+mg(a,b)\\
&=& g^{(1)}(\UU a,b)-g^{(1)}(a,\UU b)
\quad\textup{for }a,b\in\OO(K).\nonumber
\end{eqnarray}
\end{definition}

\begin{lemma}\label{t4.5}
(a) A $(T)$-structure $(H\to\C\times M,\nabla)$ induces a Higgs bundle
$(K\to M,C)$ with $K:=H|_{\{0\}\times M}$ and $C:=[z\nabla]$.

\medskip
(b) Given a $(TE)$-structure $(H\to \C\times M,\nabla)$, 
the endomorphism $\UU :=[z\nabla_{z\paa_z}]:K\to K$ is an
Euler endomorphism on the Higgs bundle from (a).

\medskip
(c) Given a $(TP)$-structure $(H\to \C\times M,\nabla,m,P)$,
the pairing $g:=(z^{-m}\cdot P)|_{K\times K}$ on $K$ is a metric
on the Higgs bundle $(K\to M,C)$ from (a).

\medskip
(d) Given a $(TEP)$-structure $(H\to\C\times M,\nnn,m,P)$,
the data from (a), (b) and (c) form a Higgs bundle
with Euler endomorphism and metric 
$(K\to M,C,\UU ,g)$. 

\medskip
(e) A $(T)$-structure $(H\to\C\times M,\nabla)$
and a basis $\uuuu{v}=(v_1,...,v_r)$ of $H|_{V\times M}$
for some neighborhood $V\subset \C$ of $0$ induce
a good connection $D$ on the Higgs bundle $(K\to M,C)$
from (a). With the notations from Remark \ref{t4.3}
and with $\uuuu{w}:=\uuuu{v}|_{\{0\}\times M}$ 
it is given by 
\begin{eqnarray}\label{4.38}
D(\uuuu{w})&:=& \uuuu{w}\cdot \sum_{i=1}^n A_i^{(1)}{\rm d}t_i.
\end{eqnarray}

\medskip
(f) A $(TE)$-structure $(H\to\C\times M,\nabla)$
and a basis $\uuuu{v}=(v_1,...,v_r)$ of $H|_{V\times M}$
for some neighborhood $V\subset \C$ of $0$ induce
a holomorphic endomorphism $Q:K\to K$ such that
$\UU $ from (b) and $Q$ form a good pair $(\UU ,Q)$
of endomorphisms on the Higgs bundle from (a)
with good connection $D$ from (e).
With the notations from Remark \ref{t4.3}
and with $\uuuu{w}:=\uuuu{v}|_{\{0\}\times M}$, $Q$ is given by
\begin{eqnarray}\label{4.39}
Q(\uuuu{v}|_{\{0\}\times M})&:=& \uuuu{v}|_{\{0\}\times M}
\cdot (-B^{(1)}).
\end{eqnarray}

\medskip
(g) A $(TP)$-structure $(H\to\C\times M,\nabla,m,P)$
and a basis $\uuuu{v}=(v_1,...,v_r)$ of $H|_{V\times M}$
for some neighborhood $V\subset \C$ of $0$ induce
on the Higgs bundle from (a) with metric from (c)
and good connection from (e) a second order pairing
$g^{(1)}$. 
With the notations from Remark \ref{t4.3}
and with $\uuuu{w}:=\uuuu{v}|_{\{0\}\times M}$, $Q$ is given by
\begin{eqnarray}\label{4.40}
g^{(1)}(\uuuu{w}^t,\uuuu{w})&=& P^{mat,(1)}.
\end{eqnarray}

\medskip
(h) In the case of a $(TEP)$-structure $(H\to\C\times M,\nabla,m,P)$
and a basis $\uuuu{v}=(v_1,...,v_r)$ of $H|_{V\times M}$
for some neighborhood $V\subset \C$ of $0$, the data
in (a)--(g) give a Higgs bundle with good connection
and good pair of endomorphisms and metric and second
order pairing and weight $m$ $(K\to M,C,D,\UU ,Q,g,g^{(1)},m)$.
\end{lemma}

{\bf Proof:}
Throughout the proof, $M$ is supposed to be small
enough so that a coordinate system $t=(t_1,...,t_n)$
on $M$ can and will be chosen, with coordinate vector fields
$\paa_i:=\paa/\paa t_i$.
And a basis $\uuuu{v}$ of $H|_{V\times M}$
for some neighborhood $V\subset \C$ of $0$ is fixed.
The notations from Remark \ref{t4.3} are used.
The restriction of the basis $\uuuu{v}$ to 
$K=H|_{\{0\}\times M}$ is denoted by $\uuuu{w}$. 

(a) The endomorphism $C_{\paa_i}$ is given by 
$C_{\paa_i}\uuuu{w}=\uuuu{w}\cdot A_i^{(0)}$.
The part for $k=0$ of \eqref{4.17} is 
$[A_i^{(0)},A_j^{(0)}]=0$ and gives \eqref{4.30}.

(b) The endomorphism $\UU $ is given by 
$\UU (\uuuu{w})=\uuuu{w}\cdot B^{(0)}$.
The part for $k=0$ of \eqref{4.18} is 
$[A_i^{(0)},B^{(0)}]=0$ and gives \eqref{4.31}. 

(c) The metric $g$ is given by 
$g(\uuuu{w}^t,\uuuu{w})=P^{mat,(0)}$. 
The part for $k=0$ of \eqref{4.21} is 
$(A_i^{(0)})^t P^{mat,(0)}=P^{mat,(0)}A_i^{(0)}$ 
and gives \eqref{4.32}.

(d) The part for $k=0$ of \eqref{4.22} is
$(B^{(0)})^t P^{mat,(0)}=P^{mat,(0)}B^{(0)}$ 
and gives \eqref{4.33}.

(e) The connection $D$ is given by 
$D_{\paa_i}\uuuu{w}=\uuuu{w}\cdot A_i^{(1)}$. 
The part for $k=1$ of \eqref{4.17} is 
$$0=\paa_iA_j^{(0)}-\paa_jA_i^{(0)}+ [A_i^{(0)},A_j^{(1)}]
+[A_i^{(1)},A_j^{(0)}]$$ 
and gives \eqref{4.34}

(f) The endomorphism $Q$ is given by
$Q(\uuuu{w})=\uuuu{w}\cdot (-B^{(1)})$.
The part for $k=1$ of \eqref{4.18} is 
$$0= \paa_iB^{(0)}+A_i^{(0)}+
[A_i^{(0)},B^{(1)}]+[A_i^{(1)},B^{(0)}]$$
and gives \eqref{4.35}. In (b) \eqref{4.31} was proved.

(g) The pairing $g^{(1)}$ is skew-symmetric because of
\eqref{4.14}. The part for $k=1$ of \eqref{4.21} is
$$\paa_iP^{mat,(0)}
=(A_i^{(0)})^t P^{mat,(1)}-P^{mat,(1)}A_i^{(0)}
+ (A_i^{(1)})^t P^{mat,(0)}+P^{mat,(0)}A_i^{(1)}$$
and gives \eqref{4.36}.

(h) The part for $k=1$ of \eqref{4.22} is
$$mP^{mat,(0)}=(B^{(0)})^t P^{mat,(1)}-P^{mat,(1)}B^{(0)}
+  (B^{(1)})^t P^{mat,(0)}+P^{mat,(0)}B^{(1)}$$
and gives \eqref{4.37}. 
\hfill $\Box$ 

\bigskip
The next definition introduces formal $(T)$-structures
and the natural enrichments. One reason is that in
a quest for classification or normal forms, it is easier
to deal first with formal $(T)$-structures and only 
later care about the convergence in $z$. 
Another one is that there is a decomposition result
for formal $(T)$-structures over a reducible germ of
an F-manifold \cite{DH20}, but not for $(T)$-structures.

\begin{definition}\label{t4.6}
Let $M$ be a complex manifold.

\medskip
(a) The sheaf $\OO_M[[z]]$ on $M$ is defined by 
$\OO_M[[z]](V):=\OO_M(V)[[z]]$
for an open subset $V\subset M$ 
(with $\OO_M(V)$ and $\OO_M[[z]](V)$ the sections of 
$\OO_M$ and $\OO_M[[z]]$ on $V$). 
Observe that the germ $(\OO_M[[z]])_{t^0}$ for $t^0\in M$ 
consists of power series $\sum_{k\geq 0}f_kz^k$ whose 
coefficients $f_k\in\OO_{M,t^0}$ have a common 
convergence domain.
In the case of $(M,t^0)=(\C^n,0)$ we write
$\OO_{\C^n}[[z]]_0=:\C\{t,z]]$.

\medskip
(b) A {\it formal $(T)$-structure} over $M$ is a free
$\OO_M[[z]]$-module $\OO(H)$ of some finite rank $r\in\N$ 
together with a map
$\nnn$ as in \eqref{4.1} where $\OO_{\C\times M}$ is replaced
by $\OO_M[[z]]$ which satisfies properties analogous to
$\nnn$ in Definition \ref{t4.1} (a), i.e. the 
Leibniz rule for $X\in\TT_M,a\in\OO_M[[z]],s\in\OO(H)$
and the flatness condition for $X,Y\in\TT_M$.

A {\it formal $(TE)$-structure}, a {\it formal $(TP)$-structure}
and a {\it formal $(TEP)$-structure} are defined analogously:
In Definition \ref{t4.1} (b)--(d) one has to replace
$\OO_{\C\times M}$ by $\OO_M[[z]]$. 
In (c), nondegeneracy of $z^{-m}\cdot P$ is only
required at $(0,t)\in\{0\}\times M$. 
\end{definition}

\begin{remarks}\label{t4.7}
(a) The formulas in the Remarks \ref{t4.3} hold also 
for formal $(T)$-structures, formal $(TE)$-structures,
formal $(TP)$-structures and formal $(TEP)$-structures
if one replaces $\OO_{\C\times M}$, $\OO_{V\times M}$
and $\OO_{(\C,0)\times M}$ by $\OO_M[[z]]$. 

\medskip
(b) Lemma \ref{t4.5} is also valid if one starts with  
formal $(T)$-structures, formal $(TE)$-structures,
formal $(TP)$-structures and formal $(TEP)$-structures.
In Lemma \ref{t4.5} (e)--(h) $\uuuu{v}=(v_1,...,v_r)$
is a $\OO_M[[z]]$-basis of $\OO(H)$. 
\end{remarks}

\section{From $(T)$-structures to pure $(TL)$-structures}
\label{c5}
\setcounter{equation}{0}

\noindent 
Definition \ref{t4.1} (e)--(h) presented extensions
to $\{\infty\}\times M$ of the vector bundles in 
Definition \ref{t4.1} (a)--(d) with good properties
along $\{\infty\}\times M$. We discuss these extensions here.
We start with the following basic existence result
Theorem \ref{t5.1} (a) 
which allows to extend (in many ways) a $(T)$-structure
to a pure $(TL)$-structure.
Theorem \ref{t5.1} (b) gives a version for $(TP)$-structures.
Theorem \ref{t5.1} (c) and (d) are weaker.
They start already with pure structures over
$\{t^0\}$ and give just extensions to $(M,t^0)$. 
A $(TE)$-structure over a point is not always
extendable to a pure $(TLE)$-structure over a point.
This is discussed in the Remarks \ref{t5.4} and Example 
\ref{t5.5}.

\begin{theorem}\label{t5.1}
(a) Let $(H\to\C\times M,\nabla)$ be a (T)-structure with 
$r=\rk H\in\N$. Let $t^0\in M$. 
Let $\uuuu{v}^0=(v_1^0,...,v_r^0)$
be a $\C\{z\}$-basis of the germ at 0 of the sheaf 
$\OO(H|_{\C\times\{t^0\}})$ of sections of the restriction 
of $H$ to $\C\times\{t^0\}$.

Then for a suitable neighborhood $U \subset M$ of $t^0$ 
an extension $(\widehat{H}\to \P^1\times U,\nabla)$ 
of $(H|_{\C\times U},\nabla)$ exists which is a pure 
(TL)-structure and such that $\uuuu{v}^0$ extends 
to a basis of global sections of $\widehat{H}$.

\medskip
(b) If the $(T)$-structure in (a) is enriched by a pairing
$P$ to a $(TP)$-structure of some weight $m\in\Z$,
one can choose the basis in $\uuuu{v}^0$ such that  
$P^{mat}(t^0)=P^{mat,(0)}(t^0)$. 
The extension in (a) induced by this basis is a pure
$(TLP)$-structure. 

\medskip
(c) Let us start in Theorem \ref{t5.1} (a) 
with a $(TE)$-structure and a global basis $\uuuu{v}^0$ of 
$H|_{\C\times \{t^0\}}$
which defines an extension to a bundle on $\P^1$
with a logarithmic pole at $\infty$, 
so a pure $(TLE)$-structure over the point $t^0$.
Then the bundle $(\widehat{H},\nnn)$ in Theorem 
\ref{t5.1} (a) is a pure $(TLE)$-structure. 

\medskip
(d) Let us start in Theorem \ref{t5.1} (a) 
with a $(TEP)$-structure and a global basis $\uuuu{v}^0$ of 
$H|_{\C\times \{t^0\}}$ with the properties in (b) and (c).
Then the bundle $(\widehat{H},\nnn)$ in Theorem 
\ref{t5.1} (a) is together with $P$ a pure $(TLEP)$-structure.
\end{theorem}

The proof will come after Corollary \ref{t5.3}.
In fact, the parts (c) and (d) are known. They can 
be derived from \cite[VI Theorem 2.1 and Proposition 2.7]{Sa02}.
But the short proofs below of (c) and (d) are convenient.
The Remarks \ref{t5.2} tell how the formulas in 
Remark \ref{t4.3} specialize in the case of pure 
$(TL)$-structures and enrichments of them.

\begin{remarks}\label{t5.2}
(i) We start with the richest structure and then
say which pieces have to be dropped in the cases
of weaker structures. 

Let $(H\to\P^1\times M,\nnn,m,P)$ be a pure
$(TLEP)$-structure of weight $m\in\Z$. Choose $M$
sufficiently small. Then a 
global basis $\uuuu{v}$ of sections on $H$ on $\P^1\times M$
exists whose restrictions to 
${H}|_{\{\infty\}\times M}$ are
$\nabla^{res,\infty}$-flat. Here
$\nabla^{res,\infty}$ denotes the restriction 
to $H|_{\{\infty\}\times M}$ of the connection of the
pure $(TL)$-structure. 
The connection $\nabla$ takes with respect to this basis 
the shape 
\begin{eqnarray}
\nabla\uuuu{v} &=& \uuuu{v}\cdot \Omega,\nonumber \\
\Omega &=& \sum_{i=1}^r z^{-1}\cdot A_i^{(0)}(t){\rm d}t_i
+ (z^{-2}B^{(0)}+z^{-1}B^{(1)}){\rm d}z.
\label{5.1}
\end{eqnarray}
Here $A_i=A_i^{(0)}$ because the restriction of the
basis $\uuuu{v}$ to $H|_{\{\infty\}\times M}$ is
$\nnn^{res,\infty}$-flat.
And $B=B^{(0)}+zB^{(1)}$ because the pole along 
$\{\infty\}\times M$ is logarithmic. Therefore the 
flatness conditions $0=d\Omega+\Omega\land\Omega$
in \eqref{4.17} and \eqref{4.18} boil here down to
\begin{eqnarray}
[A_i^{(0)},A_j^{(0)}]&=&0,\label{5.2}\\
\paa_i A_j^{(0)}-\paa_jA_i^{(0)}&=&0,\label{5.3}\\
{}[A_i^{(0)},B^{(0)}]&=&0,\label{5.4}\\
\paa_i B^{(0)}+A_i^{(0)}+[A_i^{(0)},B^{(1)}]&=&0,\label{5.5}\\
\paa_i B^{(1)}&=& 0.\label{5.6}
\end{eqnarray}
The pairing takes the shape
\begin{eqnarray}\label{5.7}
z^{-m}\cdot P(\uuuu{v}^t,\uuuu{v}) &=& P^{mat,(0)}.
\end{eqnarray}
Here $P^{mat}=P^{mat,(0)}$, because the basis
is global and the matrix of values $P^{mat}$ is holomorphic
on $\P^1\times M$, so constant in the direction of $z$.
The flatness conditions in \eqref{4.21} and \eqref{4.22}
for $P$ boil here down to 
\begin{eqnarray}\label{5.8}
(A_i^{(0)})^tP^{mat,(0)}-P^{mat,(0)}A_i^{(0)}&=& 0,\\
\paa_i P^{mat,(0)}&=& 0,\label{5.9}\\
(B^{(0)})^tP^{mat,(0)}-P^{mat,(0)}B^{(0)}&=& 0,\label{5.10}\\
(-m)P^{mat,(0)}+(B^{(1)})^tP^{mat,(0)}+P^{mat,(0)}B^{(1)}&=&0.
\label{5.11}
\end{eqnarray}

(ii) In the case of a pure $(TLE)$-structure, one drops
$P^{mat}$ and $g$ and $g^{(1)}$.
In the case of a pure $(TLP)$-structure, one drops
$B$ and $\UU $ and $Q$.
In the case of a pure $(T)$-structure, one drops
$P^{mat}$ and $B$ and $g$, $g^{(1)}$, $\UU $ and $Q$. 
\end{remarks}

The formulas in Remark \ref{t5.2} (i) tell immediately
how the data in Lemma \ref{t4.5} specialize in the
case of a pure $(TL)$-structure and enrichments of it.

\begin{corollary}\label{t5.3}
Let $(H\to\P^1\times M,\nnn,m,P)$ be a pure
$(TLEP)$-structure of weight $m\in\Z$. Choose $M$
sufficiently small. Consider a basis $\uuuu{v}$
as in Remark \ref{t5.2} (i), and 
consider the data $D,Q,g^{(1)}$ in Lemma \ref{t4.5} (e)--(h)
which are induced by this basis $\uuuu{v}$.

(a) They do not depend on this basis.

(b) The connection $D$ in Lemma \ref{t4.5} (e) is  flat, 
and the basis $\uuuu{w}=\uuuu{v}|_{\{0\}\times M}$
of the Higgs bundle $K=H|_{\{0\}\times M}$ is $D$-flat,
$D(\uuuu{w})=0$. 
The skew-symmetric second order pairing $g^{(1)}$
is equal to 0. The metric $g$ is $D$-flat, $D(g)=0$.
The endomorphism $Q+\frac{m}{2}\id$ is $D$-flat, 
$D(Q)=0$, and antisymmetric,
\begin{eqnarray}\label{5.12}
g((Q+\frac{m}{2}\id)a,b)= -g(a,(Q+\frac{m}{2}\id)b)
\quad \textup{for }a,b\in \OO(H).
\end{eqnarray}

(c) In the case of a pure $(TLE)$-structure, one drops
$g$ and $g^{(1)}$.
In the case of a pure $(TLP)$-structure, one drops
$\UU $ and $Q$.
In the case of a pure $(T)$-structure, one drops
$g$, $g^{(1)}$, $\UU $ and $Q$.
\end{corollary}

{\bf Proof:} (a) Any other basis with the same
properties consists of linear combinations of $\uuuu{v}$.
This does not change the data in Lemma \ref{t4.5}.

(b) This follows from the formulas in Remark \ref{t5.2} (i).
\hfill$\Box$ 

\bigskip

{\bf Proof of Theorem \ref{t5.1}:} 
(a) The basis $\uuuu{v}^0$ of the germ at 0 of 
$H|_{\C\times\{t^0\}}$ is a basis of $H|_{U_1\times\{t^0\}}$ 
for a suitable neighborhood $U_1\subset\C$ of $0$. 
Let $U_2\subset M$ be a contractible neighborhood of $t$.
We consider two different extensions of the sections $v^0_i$.

The sections $v^0_i$ extend uniquely to sections
$v^{flat}_i$ of $H|_{(U_1-\{0\})\times U_2}$
with $\nabla v^{flat}_i=0$. Together they form a basis
$\uuuu{v}^{flat}$ of $H|_{(U_1-\{0\})\times U_2}$.

On the other hand, holomorphic extensions $v^{hol}_i$ exist
which together form a basis of the germ of $H$ at $(0,t^0)$.

Suppose that $U_1$ and $U_2$ are chosen small enough so that
$\uuuu{v}^{hol}=(v^{hol}_1,..,v^{hol}_r)$
is a basis of $H|_{U_1\times U_2}$. Consider the base change matrix
$\Psi\in GL_r(\OO_{(U_1-\{0\})\times U_2})$ with
\begin{eqnarray}\label{5.13}
\uuuu{v}^{hol}&=& \uuuu{v}^{flat}\cdot\Psi
\end{eqnarray}  
It satisfies $\Psi(z,t^0)={\bf 1}_r$. 
By Birkhoff factorization (e.g. \cite[Proposition 4.1]{Ma83b}) 
and possibly after shrinking $U_2$, 
there are unique matrices 
$\Psi^0\in GL_r(\OO_{U_1\times U_2})$ and 
$\Psi^\infty \in GL_r(\OO_{(\P^1-\{0\})\times U_2})$ such that
\begin{eqnarray}\label{5.14}
\Psi=\Psi^\infty\cdot (\Psi^0)^{-1},\quad \Psi^\infty(\infty,t)={\bf 1}_r
\textup{ for any }t\in U_2.
\end{eqnarray}
Define the tuple of sections  
\begin{eqnarray}\label{5.15}
\uuuu{v}&:=& \uuuu{v}^{flat}\cdot \Psi^\infty
=\uuuu{v}^{hol}\cdot \Psi^0.
\end{eqnarray}
It is a holomorphic basis of $H|_{U_1\times U_2}$ (just as $\uuuu{v}^{hol}$).

We claim that the matrix valued connection 1-form $\Omega$ with 
$\nnn\uuuu{v}=\uuuu{v}\cdot\Omega$ satisfies \eqref{5.1}
(without $B^{(0)}$ and $B^{(1)}$). 
The reason is 
\begin{eqnarray*}
\nabla \uuuu{v} &=&\uuuu{v}^{flat}\cdot d\Psi^\infty
=\uuuu{v}\cdot (\Psi^\infty)^{-1}\cdot d\Psi^\infty\\
&=& \nabla(\uuuu{v}^{hol})\cdot \Psi^0 + \uuuu{v}\cdot
(\Psi^0)^{-1}\cdot d\Psi^0.
\end{eqnarray*}
The second line says that $\Omega$ has at most a pole of order 1 along
$\{0\}\times M$. The first line  says that $\Omega$ is holomorphic
on $(\P^1-\{0\})\times M$ and vanishes along $\{\infty\}\times M$.

Therefore the tuple $\uuuu{v}$ defines an extension of
$H|_{U_1\times U_2}$ to a pure (TL)-structure $\widehat{H}\to \P^1\times U_2$. 

\medskip
(b) As $P^{mat,(0)}(t^0)$ is symmetric and nondegenerate, 
the existence of a basis $\uuuu{v}^0$ with 
$P^{mat}(t^0)=P^{mat,(0)}(t^0)$ is elementary linear algebra.
 
\eqref{4.21} and $A_i=A_i^{(0)}$ and 
$P^{mat,(k)}(t^0)=0$ for $k\geq 1$ show by an induction
in the degrees of the monomials in $(t_i-t_i^0)$
in the power series expansions of the entries of
$P^{mat,(k)}$ that $P^{mat}$ is constant,
$P^{mat}=P^{mat,(0)}(t^0)$.  Therefore $z^{-m}\cdot P$
is nondegenerate at each $(z,t)\in\P^1\times M$. 
Therefore $P$ enriches the pure $(TL)$-structure to
a pure $(TLP)$-structure.

\medskip
(c) By hypothesis 
$B(z,t^0)=B^{(0)}(t^0)+zB^{(1)}(t^0)$. 
\eqref{4.18} and $A_i=A_i^{(0)}$ 
and an induction in the degrees of the monomials in 
$(t_i-t_i^0)$ in the power series expansions 
of the entries of $B^{(k)}$ for $k\geq 1$ show
$B^{(1)}=B^{(1)}(t^0)$ and $B^{(k)}=0$ for $k\geq 2$. 

\medskip
(d) This follows from (b) and (c). \hfill$\Box$

\bigskip
It is very pleasant that an extension of
a $(T)$-structure to a pure $(TL)$-structure
comes for free, and also an extension of a
$(TP)$-structure to a pure $(TLP)$-structure.
But for $(TE)$-structures and $(TEP)$-structures, 
the situation is different. 
Theorem \ref{t5.1} (c) (and (d)) gives only a relative
result: a $(TE)$-structure whose restriction to $t^0$ 
extends to a pure $(TLE)$-structure
extends also over $(M,t^0)$ to a pure $(TLE)$-structure.

\begin{remarks}\label{t5.4} (Birkhoff problem) 
(i) But there are $(TE)$-structures on $M=\{t^0\}$
which do not allow an extension to a pure $(TLE)$-structure.
The extension problem is a special case of the Birkhoff
problem, which itself is a special case of
the Riemann-Hilbert-Birkhoff problem.
The book \cite{AB94} and chapter IV in the book \cite{Sa02}
are devoted to these problems and results on them.

\medskip
(ii) An important solution of the problem of extending
a $(TE)$-structure over $M=\{t^0\}$ to a pure
$(TLE)$-structure was given implicitly 
by M. Saito in \cite{SaM89}.
Explicitly Saito solved in \cite{SaM89} 
a Riemann-Hilbert-Birkhoff problem for a 
Fourier-Laplace dual of a $(TE)$-structure with a
regular singular pole at 0. This solution is also
explained in \cite[7.4]{He02}.
Sabbah \cite{Sa98}\cite[IV 5.b.]{Sa02} 
gave a similar result for an arbitrary $(TE)$-structure. 

Both results give certain sufficient conditions
for the extendability of a $(TE)$-structure over $M=\{t^0\}$
to a pure $(TLE)$-structure.
The conditions are formulated in terms of existence of a 
second filtration which is monodromy invariant and opposite 
to a first filtration which is defined by the $(TE)$-structure. 
Both filtrations live on the space of global flat multivalued
sections of $H|_{\C^*}$. 

In Saito's result, the first filtration 
arises from comparing the regular singular germ 
$\OO(H)_0$ with a V-filtration at 0.
In Sabbah's result, the first filtration arises from
comparing the space $\OO(H)(\C)$ of global sections 
with a V-filtration at $\infty$. 

In the geometric cases (isolated hypersurface
singularities respectively tame functions on affine 
manifolds) the first filtrations are Hodge filtrations 
of mixed Hodge structures. They have indeed opposite
filtrations which are monodromy invariant. 

If the $(TE)$-structure has a logarithmic pole at 0, 
the filtrations in Saito's and Sabbah's result
coincide and the existence of an opposite and monodromy
invariant filtration is a sufficient and necessary condition
for the extendability of the $(TE)$-structure to a pure 
$(TLE)$-structure. 

\medskip
(iii) Example \ref{t5.5} is the simplest family of examples of 
a $(TE)$-structure over $M=\{t^0\}$ with a logarithmic pole
at 0 where this condition is not satisfied and 
which therefore does not allow an extension to a pure 
$(TLE)$-structure. 
\end{remarks}

\begin{example}\label{t5.5}
Let $(H'\to\C^*,\nnn)$ be a holomorphic vector bundle
of rank $r=2$ with a flat connection whose monodromy
has only one eigenvalue $\lambda\in \C^*$ 
and a $2\times 2$ Jordan block.
Write the monodromy as $M^{mon}=M_{s}\cdot M_u$
with semisimple part $M_{s}$ and unipotent part $M_u$
and nilpotent part $N=\log M_u$. 
Let $A_1,A_2$ be a basis of the space of 
global flat multivalued sections with 
$N(A_2)=-2\pi i \cdot A_1$ and $N(A_1)=0$. 
Choose $\alpha\in\C$ with $e^{-2\pi i\alpha}=\lambda$
and choose $k\in\N$.  The sections
\begin{eqnarray}\label{5.16}
v_1&:=& (z\mapsto z^\alpha A_1)\quad\textup{and}\\
v_2&:=& (z\mapsto z^{\alpha+k} e^{-\log z\frac{N}{2\pi i}}A_2)
\label{5.17}
\end{eqnarray}
are univalued and holomorphic and of moderate growth. 
The sheaf $\OO(H):=\OO_\C\cdot v_1+ \OO_\C\cdot v_2$ is the
sheaf of sections of an extension of $H'$ to a 
$(TE)$-structure with a logarithmic pole at 0, 
\begin{eqnarray}\label{5.18}
\nnn_{z\paa_z}(v_1,v_2) &=& 
(v_1,v_2)\cdot \begin{pmatrix}\alpha & z^k\\ 0&\alpha+k
\end{pmatrix}
\end{eqnarray}
The filtration $\{0\}\subset \C\cdot A_1\subset
\C\cdot A_1\oplus \C\cdot A_2$ on the space of global
flat multivalued sections does not have an opposite and
monodromy invariant filtration.
Therefore the $(TE)$-structure does not allow an extension
to a pure $(TLE)$-structure. 
\end{example}

We conclude this section with a result which concerns
formal isomorphisms between holomorphic 
$(T)$-structures over an arbitrary germ $(M,t^0)$. 

\begin{theorem}\label{t5.6} 
Let $(H\to\C\times M,\nabla)$ and 
$(\www H\to\C\times M,\www \nabla)$ be two holomorphic 
(T)-structures over the same manifold $M$. 
Let $\uuuu{v}$ and $\uuuu{\www v}$ be bases
of $\OO(H)_{(0,t^0)}$ and $\OO(\www H)_{(0,t^0)}$
at the same point $t^0\in M$. Let $T\in GL_r(\OO_{M,0}[[z]])$
define a formal isomorphism between the 
(formalized with respect to $z$) germs at $(0,t^0)$ of the 
first and the second (T)-structure which maps 
$\uuuu{v}\cdot T$ to $\uuuu{\www v}$
and which satisfies $T|_{t=t^0}\in GL_r(\C\{z\})$. 

Then $T$ is holomorphic, 
so in $GL_r(\OO_{(\C\times M,(0,t^0))})$, 
and the two (T)-structures are holomorphically isomorphic.
\end{theorem}

{\bf Proof:}
We apply Theorem \ref{t5.1} (a). 
By changing the bases holomorphically, we can suppose the following:
$\uuuu{v}$ is a global basis of a pure (TL)-structure
which extends $H^{(1)}$ and which is 
$\nnn^{res,\infty}$-flat along $\{\infty\}\times M$.
Similarly, $\uuuu{\www v}$ is a global basis of a pure 
(TL)-structure which extends $\www H$
and which is ${\www \nabla}^{res,\infty}$-flat along 
$\{\infty\}\times M$.
Finally, $T|_{t=t^0}={\bf 1}_r$. Then by
Remark \ref{t5.2}
\begin{eqnarray}\label{5.19} 
A_i&=& A_i^{(0)},\quad 
{\www A}_i = {\www A}_i^{(0)}
\end{eqnarray}
and by \eqref{4.25}
\begin{eqnarray}\label{5.20}
0&=& z\paa_i T + A_i^{(0)}T-T{\www A}_i^{(0)}.
\end{eqnarray}

An inductive argument will show $T={\bf 1}_r$. 
Suppose $t^0=0$. 
For $\uuuu{a}= (a_1,...,a_n)\in\N_0^n$ define 
$\uuuu{t}^{\uuuu{a}}:= t_1^{a_1}...t_n^{a_n}$ and 
$|\uuuu{a}|:=\sum_{i=1}^n a_i\in\N_0$.
Write
\begin{eqnarray*}
T&=&\sum_{(k,\uuuu{a})}T_{(k,\uuuu{a})}\cdot 
z^k\uuuu{t}^{\uuuu{a}},\\
A_i^{(0)}&=& \sum_{\uuuu{a}} A_{i,(0,\uuuu{a})}\cdot 
\uuuu{t}^{\uuuu{a}},\\
\www A_i^{(0)}&=& \sum_{\uuuu{a}}\www A_{i,(0,\uuuu{a})}\cdot
\uuuu{t}^{\uuuu{a}}.
\end{eqnarray*}
The condition $T|_{t=t^0}={\bf 1}_r$ is 
\begin{eqnarray}\label{5.21}
T_{(0,\uuuu{0})}={\bf 1}_r,\quad T_{(k+1,\uuuu{0})}=0
\quad \textup{for }k\geq 0.
\end{eqnarray}

The coefficient of $z^{k+1}\uuuu{t}^{\uuuu{a}}$ in \eqref{5.20}
is (with $e_i=(\delta_{ij})_{j=1,...,n}\in\N_0^n$) 
\begin{eqnarray}\label{5.22}
0&=& (a_i+1)T_{(k,\uuuu{a}+e_i)} \\
&& + \sum_{\uuuu{b}:\, 0\leq b_j\leq a_j
\, \forall\, j}A_{i,(0,\uuuu{b})}\cdot T_{(k+1,\uuuu{a}-\uuuu{b})} \nonumber\\
&& - \sum_{\uuuu{b}:\, 0\leq b_j\leq a_j\, \forall\, j} 
T_{(k+1,\uuuu{a}-\uuuu{b})}
\www A_{i,(0,\uuuu{b})}. \nonumber
\end{eqnarray}

Now we define an incomplete order $\prec$ on the set 
$\{(k,\uuuu{a})\,|\, k\in\N_0,\uuuu{a}\in\N_0^n\}$
as follows: 
\begin{eqnarray}\label{5.23}
(k,\uuuu{a})\prec(l,\uuuu{b})&\iff& 
\left\{\begin{array}{l}k+|a|<l+|b|\textup{ or}\\
k+|a|=l+|b|\textup{ and }k>l.\end{array}\right. 
\end{eqnarray}
We say that $T_{(k,\uuuu{a})}$ is {\it older} than
$T_{(l,\uuuu{b})}$ if $(k,\uuuu{a})\prec (l,\uuuu{b})$. 
\eqref{5.22} implies that $T_{(k,\uuuu{a}+e_i)}$ 
is 0 if all older $T_{(l,\uuuu{b})}$ are 0
except $T_{(0,\uuuu{0})}$ (which is ${\bf 1}_r$). 
So, inductively $T={\bf 1}_r$ follows.
\hfill$\Box$ 

\bigskip
We do not know whether one can drop in Theorem
\ref{t5.6} the condition $T|_{t=t^0}\in GL_r(\C\{z\})$. 

The following corollary allows in many cases to lift
a formal classification without work to a holomorphic
classification.

\begin{corollary}\label{t5.7}
Let $(H\to\C\times M,\nabla)$ and 
$(\www H\to\C\times M,\www \nabla)$ be two holomorphic 
(TE)-structures over the same manifold $M$. 
Let $t^0\in M$. Suppose that the restriction of 
at least one of the two $(TE)$-structures to $t^0$ is 
regular singular.

Then any formal isomorphism between the $(TE)$-structures
over the germ $(M,t^0)$ is holomorphic. 
\end{corollary}

{\bf Proof:}
Choose bases $\uuuu{v}$ and $\uuuu{\www v}$ of
$\OO(H)_{(0,t^0)}$ and $\OO(\www H)_{(0,t^0)}$.
At least one of them has moderate growth when
restricted to $(\C,0)\times\{t^0\}$.
Then the matrix $T$ of a formal isomorphism
which maps $\uuuu{v}\cdot T$ to $\uuuu{\www v}$
satisfies $T|_{t=t^0}\in GL_r(\C\{z\})$.
One applies Theorem \ref{t5.6}.
\hfill$\Box$

\section{Freedom and constraints in the steps from 
F-manifolds to Frobenius manifolds}\label{c6}
\setcounter{equation}{0}

Theorem \ref{t6.6} below  recalls a well-known recipe
for the construction of flat F-manifolds without/with
Euler field and without/with metric. 
In the case of a flat F-manifold, one builds up structures 
in 4 steps: 
\begin{list}{}{}
\item[(I)(a)] The F-manifold $M$.
\item[(II)(a)] A $(T)$-structure over the F-manifold $M$ 
(Definition \ref{t6.5}). 
\item[(III)(a)] An extension of the $(T)$-structure to a pure $(TL)$-structure.
\item[(IV)(a)] A choice of a {\it primitive section} which allows
to shift $\nnn^{res,\infty}$  to $TM$. 
This gives the flat structure of the flat F-manifold.
\end{list}
This recipe had been used in the singularity case
\cite{SaM89} (see also \cite{He02}\cite{Sa02}) 
for the construction of Frobenius manifolds with Euler fields.

The recipe has the same 4 steps in all 4 cases,
flat F-manifolds without/with Euler fields and without/with
metrics. But the freedom and the constraints in the single steps
are quite different for the 4 cases. 
One has the most freedom in the case of flat F-manifolds
without Euler field and without metric,  and the most
constraints in the case of Frobenius manifolds with
Euler fields. 
The Remarks \ref{t6.11} discuss how much freedom and 
constraints one has in the 4 steps in each of the 4 cases.

Before the recipe, Lemma \ref{t6.2} and Corollary \ref{t6.4} 
connect Higgs bundles (and enrichments) and $(T)$-structures
(and enrichments) with F-manifolds (and enrichments).

\begin{definition}\label{t6.1}
A Higgs bundle $(K\to M,C)$ is {\it primitive} 
if a section $\zeta\in\OO(K)$ exists such that the map
\begin{eqnarray}
C_\bullet\zeta:\TT_M\to \OO(K), \quad X\mapsto C_X\zeta,\label{6.1}
\end{eqnarray}
is an isomorphism. 
Such a section is called {\it primitive}.
Especially, then $\rk K=\dim M$.
\end{definition}

\begin{lemma}\label{t6.2}
(a) \cite[Lemma 4.1]{He03}
A primitive Higgs bundle $(K\to M,C)$ 
induces a unique multiplication $\circ$
on $TM$ with $C_{X\circ Y}=C_XC_Y$. With respect to an arbitrary
primitive section $\zeta$, it is uniquely determined by 
$C_{X\circ Y}\zeta =C_XC_Y\zeta$. Also a unique unit field 
$e\in\TT_M$ with $e\circ =\id$ is induced. 
It is determined by $C_e\zeta=\zeta$.

\medskip
(b) \cite[Lemma 4.3]{He03}
A primitive Higgs bundle with good connection 
$(K\to M,C,D)$ induces an F-manifold structure on $M$.

\medskip
(c) \cite[Lemma 4.3]{He03}
A primitive Higgs bundle with good connection and a good pair 
of endomorphisms $(K\to M,C,D,\UU ,Q)$ induce an Euler field $E$ 
on the F-manifold $M$ (from (b)) with $C_E=-\UU $.
With respect to an arbitrary primitive section $\zeta$, 
the Euler field $E$ is determined by $C_E\zeta=-\UU \zeta$.
\end{lemma}

\begin{remarks}\label{t6.3}
(i) Let $(K\to M,C)$ be a primitive Higgs bundle.
For each $t^0\in M$, vectors in an open set in $T_{t^0}M$
extend to primitive sections on $K|_{(M,t^0)}$.
So, primitive sections are not at all unique.

(ii) Any choice of a primitive section $\zeta$ yields an
isomorphism \eqref{6.1}. All enrichments in Definition 
\ref{t4.4} can be pulled back from $K$ to $TM$ with this 
isomorphism and satisfy there the same properties. 
Let us call them $C^{M}$, $\UU^{M}$, 
$g^{M,\zeta}$, $D^{M,\zeta}$, $Q^{M,\zeta}$, 
$g^{(1),M,\zeta}$. They will be used for a special choice
of $\zeta$ in Theorem \ref{t6.6}.
Here $C^{M}$ and $\UU^{M}$ do not depend
on the choice of $\zeta$, all others depend on the choice
of $\zeta$. And here $C^{M}_X=X\circ$ and 
$-\UU^{M}=E\circ$. 
\end{remarks}

Lemma \ref{t6.2} and Lemma \ref{t4.5} have the following
corollary. 

\begin{corollary}\label{t6.4}
(a) Let $(H\to\C\times M,\nnn)$ be $(T)$-structure
such that the underlying Higgs bundle 
$(K=H|_{\{0\}\times M}\to M,C)$ is primitive. 
Then $M$ becomes an F-manifold.

(b) Let $(H\to\C\times M,\nnn)$ be a $(TE)$-structure
such that the underlying Higgs bundle is primitive.
Then $M$ becomes an F-manifold with Euler field $E$. 
\end{corollary}

\begin{definition}\label{t6.5}
A {\it $(T)$-structure} (or any enrichment in Definition 
\ref{t4.1}) {\it over an F-manifold} means a $(T)$-structure
such that the underlying Higgs bundle is primitive,
and then the induced F-manifold structure on the base
space is considered. In the case of a $(TE)$-structure
also the induced Euler field is considered.
\end{definition}

The following theorem is a recipe for constructing
flat F-manifolds without/with Euler field and 
without/with metric. 
It was first used essentially by M. Saito
for the construction of Frobenius manifold structures 
on base spaces of universal unfoldings of 
isolated hypersurface singularities \cite{SaM89}.
The version for Frobenius manifolds with Euler field 
is given in \cite[VII Theorem 3.6]{Sa02}
and in \cite[Theorem 5.12]{He03}.
We recall the proof. The proof is short because of the 
preparations in the sections \ref{c4} and \ref{c5}.

\begin{theorem}\label{t6.6}
(a) Let $(H\to\P^1\times M,\nabla)$ be a pure (TL)-structure
over an F-manifold. 
Let $\omega\in \OO(H)(\P^1\times M)$ be a global section 
such that its restriction $\omega|_{\{\infty\}\times M}$ 
to $H|_{\{\infty\}\times M}$ is flat with respect to the restriction 
$\nabla^{res,\infty}$ of 
$\nabla$ to $H|_{\{\infty\}\times M}$ and such that its restriction 
$\zeta:=\omega|_{\{0\}\times M}$ to $K$
is a primitive section of the primitive Higgs bundle.
Consider the connection $D$ on $K$ from Corollary \ref{t5.3}.
The isomorphism
\begin{eqnarray}\label{6.2}
(C_\bullet\zeta)^{-1}:K\to TM
\end{eqnarray}
maps $D$ to a connection $D^{M,\zeta}$ on $TM$ such that
$(M,\circ,e,D^{M,\zeta})$ is a flat F-manifold.

\medskip
(b) Let $(H\to\P^1\times M,\nnn)$ be a pure $(TLE)$-structure
over an F-manifold. By Corollary \ref{t6.5} the F-manifold
$M$ has an Euler field $E$. 
Consider the endomorphism $Q$ from Corollary \ref{t5.3}.
Let $\omega$ and $\zeta$ be as in (a) and suppose that $\zeta$
is a section of eigenvectors of $Q$,
\begin{eqnarray}\label{6.3}
Q(\zeta)&=& \frac{d}{2}\cdot \zeta \quad\textup{for some }d\in\C.
\end{eqnarray}
Then $(M,\circ,e,E,D^{M,\zeta})$ is a flat F-manifold 
with Euler field.

\medskip
(c) Let $(H\to\P^1\times M,\nnn,m,P)$ be a pure $(TLP)$-structure
of weight $m\in\Z$ over an F-manifold. 
Let $\omega$ and $\zeta$ be as in (a).
Shift the metric $g$ from Lemma \ref{t4.5} (c)
with the isomorphism in \eqref{6.2} to a metric $g^{M,\zeta}$ 
on $TM$. 
Then $(M,\circ,e,g^{M,\zeta})$ is a Frobenius manifold, 
and $D^{M,\zeta}$ is the Levi-Civita connection of 
$g^{M,\zeta}$.

\medskip
(d) Let $(H\to\P^1\times M,\nnn,m,P)$ be a pure $(TLEP)$-structure
of weight $m\in\Z$ over an F-manifold. 
Let $\omega$ and $\zeta$ be as in (a) and suppose that $\zeta$
is a section of eigenvectors of $Q$ with eigenvalue $\frac{d}{2}$
as in \eqref{6.3}.
Then $(M,\circ,e,E,g^{M,\zeta})$ is a Frobenius manifold 
with Euler field and 
$\Lie_E(g^{M,\zeta})=(2-d-m)\cdot g^{M,\zeta}$.
\end{theorem}

{\bf Proof:}
(a) Here $C^M$ is induced by $C$ via the isomorphism in \eqref{6.2}.
Therefore the potentiality $D^{M,\zeta}(C^M)=0$ follows from
the potentiality $D(C)=0$. The flatness of $D^{M,\zeta}$ follows from the
flatness of $D$ in Corollary \ref{t5.3}. 
A local basis $\uuuu{v}$ as in Corollary \ref{t5.3} can be chosen
such that $\omega$ is part of the basis. 
Thus $D(\uuuu{v}|_{\{0\}\times M})=0$ implies especially
$D(\zeta)=0$. As $e$ is the image of $\zeta$ under the isomorphism
in \eqref{6.2}, this gives $D^{M,\zeta}(e)=0$. This and the potentiality
$0=D^{M,\zeta}(C^M)$ applied to $e$ in the explicit version 
in \eqref{3.3} give the torsion freeness of $D^{M,\zeta}$.

(b) The pair $(-E\circ, Q^{M,\zeta})$ inherits from the pair $(\UU ,Q)$
the property of being a good pair of endomorphisms. Thus
\begin{eqnarray}\label{6.4}
0=D^{M,\zeta}(-E\circ)-[C^M,Q^{M,\zeta}]+C^M.
\end{eqnarray}
The unit field $e$ inherits from $\zeta$ the property 
of being a section of eigenvectors of $Q^{M,\zeta}$, $
Q^{M,\zeta}(e)=\frac{d}{2}\cdot e$. Inserting $X$ and $e$ into 
\eqref{6.4} gives
\begin{eqnarray}
D^{M,\zeta}_X(-E)&=& D^{M,\zeta}_X(-E\circ)(e)= [X\circ,Q^{M,\zeta}](e)-X\circ e\nonumber\\
&=&X\circ Q^{M,\zeta}(e)-Q^{M,\zeta}(X\circ e)-X \nonumber\\
&=& -Q^{M,\zeta}(X)-\frac{2-d}{2}X, \nonumber\\
\textup{so } D^{M,\zeta}_\bullet(E)&=&  Q^{M,\zeta}+\frac{2-d}{2}\id .\label{6.5}
\end{eqnarray}
$Q$ in Corollary \ref{5.3} is flat, $D(Q)=0$, and 
$Q^{M,\zeta}$ inherits this property, $D^{M,\zeta}(Q^{M,\zeta})=0$.
Therefore $D^{M,\zeta}_\bullet(E)$ is flat.

(c) $g$ in Corollary \ref{t5.3} is flat, $D(g)=0$,
and $g^{M,\zeta}$ inherits this property, $D^{M,\zeta}(g^{M,\zeta})=0$.
As $D^{M,\zeta}$ is torsion free by part (a), $D^{M,\zeta}$ is the Levi-Civita
connection of $g^{M,\zeta}$. 
The metric $g^{M,\zeta}$ is multiplication invariant because of 
Lemma \ref{t4.5} (c) and Definition \ref{t4.4} (c).
This and part (a) give part (c).

(d) Everything except $\Lie_E(g^{M,\zeta})
=(2-d-m)\cdot g^{M,\zeta}$ follows from
(b) and (c). Let $X,Y\in\TT_M$ be flat vector fields.
The torsion freeness of $D^{M,\zeta}$ gives 
$[E,X]=-D_X^{M,\zeta}E$, $[E,Y]=-D_Y^{M,\zeta}E$.
\begin{eqnarray*}
&& \Lie_E(g^{M,\zeta})(X,Y)\\
&=& E(g^{M,\zeta}(X,Y))
-g^{M,\zeta}([E,X],Y)-g^{M,\zeta}(X,[E,Y])\\
&=& 0+g^{M,\zeta}(D_X^{M,\zeta}E,Y)+g^{M,\zeta}(X,D_Y^{M,\zeta}E)\\
&\stackrel{\eqref{6.5}}{=}& 
g^{M,\zeta}(Q^{M,\zeta}X +\frac{2-d}{2}X,Y)
+ g^{M,\zeta}(X,Q^{M,\zeta}Y+\frac{2-d}{2}Y)\\
&\stackrel{\eqref{5.12}}{=}& (2-d-m)g^{M,\zeta}(X,Y).
\end{eqnarray*}
\hfill$\Box$

\bigskip

Before we finish this section with a discussion of 
the freedom and the constraints in the 4 steps 
at the beginning of this section for the construction 
of flat F-manifolds and enrichments,
two other results are given.

\begin{theorem}\label{t6.7}
(a) Let $(K\to M,C,D)$ be a primitive Higgs bundle with a good 
connection. By Lemma \ref{t6.2} (b), it induces an F-manifold
structure on $M$. 
Let $\II_M\subset \OO(T^*_M)$ be the ideal which defines
the analytic spectrum of the multiplication
(see Definition \ref{t2.11}, also for the Poisson bracket
$\{.,.\}$ on $\OO(T^*_M)$). Then
\begin{eqnarray}\label{6.6}
\{\sqrt{\II_M},\sqrt{\II_M}\} &\subset & \sqrt{\II_M}.
\end{eqnarray}
(b) Let $(H\to\C\times M,\nnn)$ be a $(T)$-structure over
an F-manifold. Then \eqref{6.6} holds.
\end{theorem}

{\bf Proof:} (a) By Theorem \ref{t2.3}, there is a map 
$l:M\to\N$
such that at each point $t^0\in M$ the germ of the F-manifold
decomposes into a product $\prod_{k=1}^{l(t^0)}(M_k,t^{0,k})$
of $l(t^0)$ many irreducible germs of F-manifolds. 
The caustic $\KK$ is the set 
\begin{eqnarray}\label{6.7}
\KK&:=& \{t\in M\,|\, l(t)\textup{ is not the maximal value}\}.
\end{eqnarray}
This generalizes the notion of a caustic in the generically
semisimple case. The caustic $\KK$ is empty or a hypersurface 
\cite[Proposition 2.6]{He02}.

\eqref{6.6} is a local property. It is sufficient to prove it
near a point $t^0\in M-\KK$. Consider such a point. 
Write $l:=l(t^0)$.
Over the germ $(M,t^0)=\prod_{k=1}^l(M_k,t^{0,k})$,
also the bundle $K$ decomposes uniquely into a direct sum 
$\bigoplus_{k=1}^l K_k$ of subbundles. 
Here $K_k =C_{e_k}(K)$ is the image of the projection 
$C_{e_k}$, where $e_k\in\TT_M$ is the lift to $\TT_M$
of the unit field of $(M_k,t^{0,k})$. 

Restrict $K_k$ to 
$(M_k,t^{0,k})\times \prod_{j\neq k}\{t^{0,j}\}$
and call the resulting bundle $K^{(k)}$.
The Higgs field $C$ restricts to a Higgs field $C^{(k)}$
on this bundle, and the pair $(K^{(k)}\to M_k,C^{(k)})$
is again a primitive Higgs bundle 
(with $M_k$ a suitable representative of the germ).
Also the connection $D$ has a part 
$D^{(k)}:\OO(K^{(k)})\to \Omega^1_{M_k}\otimes\OO(K^{(k)})$
which acts on $K^{(k)}$. 
This connection inherits from $D$ the property of being good.
Also the cotangent bundle and the Poisson bracket
split according to the decomposition 
$(M,t^ 0)=\prod_{k=1}^{l}(M_k,t^{0,k})$.
Therefore we can and will restrict to the case $l=1$. 

Let $t=(t_1,...,t_n)$ be coordinates on $(M,t^0)$
with $t^0=0$ and $\paa_i:=\paa/\paa t_i$ and $e=\paa_1$.
Replace $M$ be a neighborhood of $0$ in $M-\KK$
on which the coordinates are defined. 
Then $\KK=\emptyset$, and for any $t\in M$ 
the algebra $T_tM$ is irreducible.
Therefore $\paa_i\circ$ has for each $t$
only one eigenvalue $\rho_i(t)$, i.e. 
$(\paa_i-\rho_i\cdot \paa_1)\circ$ is for $t\in M$ nilpotent. 
The functions $\rho_i:M\to\C$ are holomorphic.
The radical of the ideal $\II_M$ is 
\begin{eqnarray}\label{6.8}
\sqrt{\II_M}&=& (y_1-1,y_2-\rho_2,...,y_n-\rho_n)
\end{eqnarray}
It satisfies \eqref{6.6} if and only if for all $i\neq j$
$\paa_i\rho_j=\paa_j\rho_i$. 

But this follows from the potentiality $D(C)=0$ in the
following way. Also $C_{\paa_i}$ has at $t\in M$ 
only the eigenvalue $\rho_i(t)$. Therefore 
$\rho_i(t)=\tr(C_{\paa_i}|_t)$. 
Choose a basis $\uuuu{w}$ of $K\to M$ near 0 and
write $C_{\paa_i}\uuuu{w}=\uuuu{w}\cdot A_i^{(0)}$
and $D_{\paa_i}\uuuu{w}=\uuuu{w}\cdot A_i^{(1)}$.
The potentiality $D(C)=0$ says
\begin{eqnarray}\label{6.9}
0&=& \paa_i A_j^{(0)}-\paa_j A_i^{(0)} + 
[A_i^{(0)},A_j^{(1)}] + [A_i^{(1)},A_j^{(0)}].
\end{eqnarray}
Taking the traces gives 
$$\paa_i\rho_j=\paa_i \tr(A_j^{(0)}) =
\tr(\paa_i A_j^{(0)}) = \tr(\paa_j A_i^{(0)}) 
= \paa_j \tr(A_i^{(0)})=\paa_j\rho_i.$$

(b) This follows from part (a) and Lemma \ref{t4.5} (e).
\hfill$\Box$

\begin{remark}\label{t6.8}
Part (b) of Theorem \ref{t6.7} was known before. 
A $(T)$-structure over an F-manifold $M$ is a special case 
of an $\RR_\XX$-structure in the sense of \cite{Sa05},
where $X=M$. Proposition 1.2.5 in \cite{Sa05} implies
\eqref{6.6}. But the proof is an application of a deep result
of Gabber. It is good to have a short and elementary proof
for Theorem \ref{t6.7}.
\end{remark}

Certain $(TE)$-structures can be {\it unfolded} uniquely
to $(TE)$-structures over F-manifolds with Euler fields.
The basic result is due to Malgrange \cite{Ma86}.
The following Theorem \ref{t6.9}
cites in part (a) and (b) 
a generalization in \cite[Theorem 2.5]{HM04}
and in part (c) Malgrange's result. 
For the notion of a {\it universal unfolding} of a 
$(TE)$-structure, we refer to \cite[Definition 2.3]{HM04}.

\begin{theorem}\label{t6.9}
(a) Let $(H\to\C\times(M,t^0),\nnn)$ be a $(TE)$-structure over
a germ $(M,t^0)$ of a manifold. 
Let $(K\to (M,t^0),C)$ be the induced Higgs bundle over
$(M,t^0)$. Suppose that a vector $\zeta_{t^0}\in K_{t^0}$ 
with the following properties exists.
\begin{list}{}{}
\item[(IC)] (Injectivity condition) The map 
$C_\bullet \zeta_{t^0}:T_{t^0}M\to K_{t^0}$ is injective.
\item[(GC)] (Generation condition) $\zeta_{t^0}$ and
its images under iteration of the maps 
$\UU|_{t^0}:K_{t^0}\to K_{t^0}$ and $C_X:K_{t^0}\to K_{t^0}$
for $X\in T_{t^0}M$ generate $K_{t^0}$.
\end{list}
Then a universal unfolding of the $(TE)$-structure over
a germ $(M\times \C^l,(t^0,0))$ ($l\in\N_0$ suitable)
exists. 
It is unique up to isomorphism. The underlying Higgs bundle
is primitive. Thus the germ $(M\times \C^l,(t^0,0))$
is an F-manifold with Euler field.

\medskip
(b) A $(TE)$-structure over the germ of an F-manifold
(automatically with Euler field) is its own universal
unfolding.

\medskip
(c) \cite{Ma86} A $(TE)$-structure over a point $t^0$ has
a universal unfolding if the endomorphism 
$[z^2\nnn_{\paa_z}]=\UU :K_{t^0}\to K_{t^0}$ 
is regular, i.e. it has for each eigenvalue only one 
Jordan block. 
In that case, the F-manifold with Euler field 
which underlies the universal unfolding, 
is by definition (Definition \ref{t2.15}) regular. 
\end{theorem}

Theorem 1.3 in \cite{DH17} gives the structure of
regular F-manifolds. 

\begin{remarks}\label{t6.10}
(i) Any $(TE)$-structure over the germ $(M,t^0)$ of a regular 
F-manifold is by Malgrange's result determined by its
restriction to $t^0$. That restriction can be treated
and understood as in \cite{Ma83a}. 
By the Hukuhara-Levelt-Turrittin theorem, a branched covering
$\C\to\C$, $z\mapsto z^k$, exists such that the pull-back 
of the underlying meromorphic bundle has a formal decomposition
into a tensor product of rank 1 bundles with irregular
connections and bundles with regular singular connections.
A Stokes structure distinguishes the meromorphic bundle
within this formal isomorphism class. 

(ii) We wonder whether an upper bound for the orders $k$ 
of the branched coverings in (i) for all $(TE)$-structures 
over $t^0$ exists, for which 
$[z^2\nnn_{\paa_z}]=\UU :K_{t^0}\to K_{t^0}$ 
is a regular endomorphism. 
If yes, the number of continuous parameters for
$(TE)$-structures over a germ of a regular F-manifold
would be bounded. 
\end{remarks}

\begin{remarks}\label{t6.11}
Here we discuss the freedom and the constraints in the
4 steps (I), (II), (III) and (IV) at the beginning of this
section for the construction of flat F-manifolds
without/with Euler fields and without/with metrics.
These cases are indexed by (a), (b), (c) and (d) as follows.

\begin{tabular}{c|c|c}
 & without Euler field & with Euler field \\ \hline
without metric & (a) flat F-manifold & (b) flat F-manifold with $E$\\
\hline
with metric & (c) Frobenius mfd. & (d) Frobenius mfd. with $E$
\end{tabular}

\medskip
(I) (a)  It is rather easy to construct F-manifolds,
see Remark \ref{t2.14}. Their structure is quite well
understood, see section \ref{c2} and \cite{He02}.

(I) (b) Existence of Euler fields for F-manifolds is 
discussed in \cite[3.2]{He02}.
Many F-manifolds do not have Euler fields. 
Some have many different Euler fields. 
There may be functional parameters.
A generically semisimple F-manifold $M$ has locally 
on the complement $M-\KK$ of the caustic 
an affine space (of dimension $\dim M$) of Euler fields, see
Remark \ref{t2.7} (i). 
But at a point $t^0\in\KK$, where the germ of the F-manifold
is irreducible, at most one of them extends to an
Euler field in a neighborhood of $t^0$, see Theorem \ref{t2.10}.

(I) (c) Given an F-manifold, locally a multiplication invariant
metric exists if and only if all algebras $T_tM$ are
Frobenius algebras, see the Remarks \ref{t2.16}.
In any F-manifold, the set of points $t\in M$ with
$T_tM$ a Frobenius algebra is empty or an open set.
Often this set is equal to $M$, but also often it is not,
see \cite[Proposition 5.32 and the Remarks 5.33]{He02}
for examples where it is not equal to $M$.

(I) (d) The examples cited in (c) are F-manifolds with
Euler fields. 

\medskip
(II) (a) If a $(T)$-structure lives over an F-manifold,
then the F-manifold satisfies \eqref{6.6} by Theorem
\ref{t6.7}. 
The F-manifolds in \cite[2.5.2 and 2.5.3]{HMT09} do not
satisfy this condition. Therefore there are no 
$(T)$-structures over them. 
On the other hand, over many F-manifolds, $(T)$-structures
live, especially over semisimple and (more generally)
regular F-manifolds. 
Conjecture \ref{t7.2} (a) says that over any irreducible 
germ of a generically semisimple F-manifold, 
$(T)$-structures live.
Though Conjecture \ref{t7.2} (a) is not precise, how many
$(T)$-structures live over such a germ.
Theorem \ref{t8.2} gives already for the $(T)$-structures
over the F-manifold $I_2(m)$ functional parameters 
(in $F_2\supset F_3$ in Remark \ref{t8.3} (i)).
And already in the case of $I_2(m)$ it is difficult to give
holomorphic normal forms for $(T)$-structures over $I_2(m)$,
see Remark \ref{t8.3} (iii) and Theorem \ref{t8.4}.

(b) $(TE)$-structures live over semisimple F-manifolds
and (more generally) regular F-manifolds.
Because of Theorem \ref{t6.9}, such a $(TE)$-structure
is determined by its restriction to a point $t^0$.
In the semisimple case, that restriction can be encoded by 
the eigenvalues of $\UU|_{t^0}$, two Stokes matrices,
and the exponents of the regular singular rank 1 pieces.
See Remark \ref{t6.10} for the case of regular F-manifolds. 

Conjecture \ref{t7.2} (b) says that over any irreducible
germ of a generically semisimple F-manifold with Euler field,
$(TE)$-structures live. And here it is precise how many
and which invariants distinguish them.
Theorem \ref{t7.4} says that Conjecture \ref{t7.2} 
is true for the 2-dimensional F-manifolds $I_2(m)$. 

On the other hand, an example of M. Saito \cite{SaM17}
can be enhanced to a family of $(TE)$-structures
over a fixed 3-dimensional everywhere irreducible
F-manifold with Euler field which is nowhere a regular
F-manifold. This family has a functional parameter
$g(t_1)\in\C\{t_1\}$. We plan to treat it in \cite{DH19-2}.

(c) and (d) We expect that the pairing does not pose
problems. If a $(T)$-structure over an F-manifold
exist, whose tangent spaces are all Frobenius algebras,
we expect that also a $(TP)$-structure exists.
Conjecture \ref{t7.2} says that over an irreducible germ of
a generically semisimple F-manifold with Euler field
exactly one $(TEP)$-structure exists if $T_{t^0}M$
is a Frobenius algebra (and no one if not).

\medskip
(III) (a) Theorem \ref{t5.1} (a) says that any $(T)$-structure
over a germ of a manifold extends in many ways
to a pure $(TL)$-structure. So, this extension comes for free. 
The freedom in Theorem \ref{t5.1} in the choice of $\uuuu{v}^0$ 
shows that one has functional parameters.

(b) On the other hand, not every $(TE)$-structure extends to
a pure $(TLE)$-structure. The problem arises already
in the case of $(TE)$-structures over a point $t^0$. 
Example \ref{t5.5} is an example of rank 2 with a 
logarithmic pole at $z=0$. 
The problem of extendability is the Birkhoff problem.
It is discussed in the Remarks \ref{t5.4}.
If the restriction of a $(TE)$-structure to a point
extends to a pure $(TLE)$-structure, the whole
$(TE)$-structure over a germ $(M,t^0)$ 
extends to a pure $(TLE)$-structure. This rather classical
fact is recalled (and reproved) in Theorem \ref{t5.1} (c).

(c) Theorem \ref{t5.1} (b) says that any $(TP)$-structure
over a germ of a manifold extends in many ways
to a pure $(TLP)$-structure. Also this extension comes for free.

(d) In the case of a $(TEP)$-structure over a point $t^0$,
the problem of extension to a pure $(TLEP)$-structure
is again a Birkhoff problem. Though the additional 
constraint by the pairing is not so serious. 
Only, it requires some extra care. 

\medskip
(IV) (a) Consider a pure $(TL)$-structure over a 
sufficiently small manifold $M$. Then the global sections
whose restrictions to $\{\infty\}\times M$
are $\nnn^{res,\infty}$-flat, form a vector space 
$\VV^{res,\infty}$ with $\dim \VV^{res,\infty}=\rk H$.
The last ingredient towards a flat F-manifold
is in Theorem \ref{t6.6} a global section 
$\omega\in\VV^{res,\infty}$ whose restriction to
$\{0\}\times M$ is a primitive section of the 
Higgs bundle. 
All sections in an open subset of $\VV^{res,\infty}$  
can take the role of $\omega$.
Therefore, one can pass for free from a pure
$(TL)$-structure on an F-manifold to a flat F-manifold,
and one has $\rk H$ parameters.

(b) In the case of a pure $(TLE)$-structure, 
$Q$ acts on $\VV^{res,\infty}$, and 
one has the additional constraint that $\omega$
must be a section of eigenvectors of $Q$.
If one has bad luck, the open subset of $\VV^{res,\infty}$
of sections whose restrictions to $\{0\}\times M$ 
are primitive sections of the Higgs bundle
does not intersect the union of eigenspaces of $Q$
in $\VV^{res,\infty}$. For example in the general
cases of tame functions on affine manifolds which are
considered in \cite{DS03}, this is not clear.
Though they intersect in the case of Newton nondegenerate
polynomials or Laurent series, and also in the case
of isolated hypersurface singularities.

(c) and (d) Including the pairing does not change
the situation. It does not give constraints on the
section $\omega$. 
\end{remarks}

\section{A conjecture on $(TE)$-structures over 
irreducible germs of generically semisimple F-manifolds
with Euler fields}\label{c7}
\setcounter{equation}{0}

\noindent
Conjecture \ref{t7.2} predicts existence of
$(T)$-structures, $(TE)$-structures, $(TP)$-structures
and $(TEP)$-structures over irreducible germs of 
generically semisimple F-manifolds.
Part (a) on $(T)$-structures and $(TP)$-structures
is not precise on the size of families and their parameters.
Part (b) is precise about the family of all $(TE)$-structures
over one germ and about their parameters,
and it predicts one or zero $(TEP)$-structures of weight
$m\in\Z$. Before giving the conjecture, we define and discuss
the {\it regular singular exponents} of a $(TE)$-structure 
over a generically semisimple F-manifold with Euler field.

\begin{remarks}\label{t7.1}
Let $(H\to\C\times M,\nnn)$ be a $(TE)$-structure over a 
generically semisimple F-manifold with Euler field,
with $n=\dim M=\rk H$. 
Recall the definition of caustic $\KK$, Maxwell stratum 
$\KK_2$ and bifurcation set $\KK^{bif}=\KK\cup\KK_2$ from
Remark \ref{t2.7} (iii).

\medskip
(i) The Hukuhara-Levelt-Turrittin theorem on the formal
decomposition of an irregular pole at $0\in\C$
of a meromorophic bundle with a meromorphic connection
allows an enhancement to $(TE)$-structures \cite{HS07}.

Especially, the $(TE)$-structure over $t^0\in M-\KK^{bif}$
has a formal decomposition into
a sum of rank 1 $(TE)$-structures. Each rank 1
$(TE)$-structure itself is a tensor product 
of an irregular rank 1 $(TE)$-structure with
$\nnn=d-u_i(t^0)z^{-2}{\rm d}z$ and a logarithmic 
rank 1 $(TE)$-structure with $\nnn=d+\alpha_i z^{-1}{\rm d}z$.
The values $\alpha_1,...,\alpha_n$ are called the 
{\it regular singular exponents} of the $(TE)$-structure 
over $t^0$. 
The values $u_1(t^0),...,u_n(t^0)$ are the
eigenvalues of $E\circ=-\UU^m:T_{t^0}M\to T_{t^0}M$ 
and of $-\UU :K_{t^0}\to K_{t^0}$.

The $(TE)$-structure over $t^0$ is determined by
the values $u_1(t^0),...,u_n(t^0)$, the
values $\alpha_1,...,\alpha_n$ and the Stokes structure
\cite{Ma83a}. 

\medskip
(iii) Locally near $t^0$, $u_1,...,u_n$ are canonical
coordinates of the semisimple F-manifold.
Moving locally in the F-manifold, the canonical coordinates 
vary, but the Stokes structure and the regular singular
exponents are constant. 

Without the notions of F-manifold and $(TE)$-structures,
Malgrange considered this situation in
\cite{Ma83b}.

\medskip
(iv) Consider the analytic spectrum $L_M\subset T^*_M$
of the F-manifold. The projection 
$L_M|_{M-\KK}\to M-\KK$ is a covering with $n$ sheets.
In a small open set $V\subset M-\KK$, 
each sheet is associated to one regular singular exponent
$\alpha_i$, and this value gives a constant function 
with value $\alpha_i$ on this sheet.
One obtains a global function $L_M|_{M-\KK}\to\C$
which is locally constant. It gives a function
from the finite set of topological components of
$L_M|_{M-\KK}$ to $\C$. These topological components 
of $L_M|_{M-\KK}$ are the restrictions of the
components of $L_M$ as a reduced complex space.
So, one obtains a map
\begin{eqnarray}\label{7.1}
\AAA^{rse}:\{\textup{components of }L_M\}\to\C
\end{eqnarray}
({\it rse} for {\it regular singular exponent}).
It is an invariant of the $(TE)$-structure over the
generically semisimple F-manifold.
\end{remarks}

\begin{conjecture}\label{t7.2}
(a) Let $(M,t^0)$ be an irreducible germ of a generically 
semisimple F-manifold. 

Then (T)-structures over $(M,t^0)$
exist. In general, there are functional parameters
(parameters in rings of power series).

If $T_{t^0}M$ is a Frobenius algebra, then
$(TP)$-structures exist. Again, in general, there are
functional parameters.

\medskip
(b) Let $(M,t^0)$ be an irreducible germ of a generically
semisimple F-manifold with Euler field, and let
$M$ be a sufficiently small representative. 
Let its analytic spectrum
$L_M=\textup{Spec}_{\OO_M}(TM)\subset T^*M$ have 
$l\in\{1,2,...,\dim M\}$
components (as a reduced complex space).

Then $(TE)$-structures over $(M,t^0)$ exist.
At (an arbitrary point) 
$t^1\in M-\KK^{bif}$ they have all the same Stokes structure.
For each map 
\begin{eqnarray}\label{7.2}
\AAA:\{\textup{components of }L_M\}\to\C,
\end{eqnarray}
precisely one $(TE)$-structure with $\AAA=\AAA^{rse}$ 
exists. Therefore the $(TE)$-structures over $(M,t^0)$  
are parametrized by $\C^l$. 

If $T_{t^0}M$ is a Frobenius algebra and $m\in\Z$, 
then the $(TE)$-structure with all regular singular exponents
equal to $m/2$ can be enriched to a $(TEP)$-structure
of weight $m$.
\end{conjecture}

\begin{remarks}\label{t7.3}
(i) Part (a) of Conjecture \ref{t7.2} is sketchy,
as we do not have a good understanding of the
parameters behind the $(T)$-structures over
one irreducible germ of a generically semisimple
F-manifold, not even in the case of the 
2-dimensional F-manifolds $I_2(m)$,
see Remark \ref{t8.3}.
Anyway, our main interest is in part (b).

\medskip
(ii) We can imagine a generalization of Conjecture
\ref{t7.2} (b) to irreducible germs of generically
regular F-manifolds. Generically regular F-manifolds 
share with generically semisimple F-manifolds 
the property that locally a $(TE)$-structure 
over the F-manifold with Euler field is determined by
its restriction to a generic single point, 
see Theorem \ref{t6.9}.

\medskip
(iii) We expect functional parameters for $(TE)$-structures
over a germ of an F-manifold with Euler field
which is nowhere regular. See Remark \ref{t6.11} (II)(b).

\medskip
(iv) We expect that the (TE)-structures in Conjecture 
\ref{t7.2} (b) have all the same Stokes structure. 
Conjecture \ref{t7.2} (b) is an existence and a 
uniqueness conjecture also for the Stokes structure.
At a point $t^0\in M-\KK^{bif}$ one can start with a semisimple
(TE)-structure on $(\C,0)\times\{t^0\}$ with arbitrary 
Stokes structure and arbitrary regular singular exponents
$\alpha_1,...,\alpha_n$. 
A local isomonodromic extension exists and is unique.
A global isomonodromic extension to $M-\KK^{bif}$ exists only if
the Stokes structure fits to itself after all possible paths
around $\KK^{bif}$ and if the $\alpha_i$ from the rank 1 
pieces over points in the same component of $L_M$ coincide.
Suppose that this holds for some Stokes structure
and some choice of $\alpha_1,...,\alpha_n$. 
Still the holomorphic extendability of the $(TE)$-structure 
from $M-\KK^{bif}$ to $\KK^{bif}$ is a nontrivial constraint. 
Conjecture \ref{t7.2} (b) claims that precisely one 
Stokes structure satisfies all constraints.

\medskip
(v) In the case of $I_2(m)$ with $m\geq 3$, we checked that 
there are $[\frac{m}{2}]$ Stokes structures which give 
$(TE)$-structures on $M-\KK^{bif}=M-\KK=\C^2-\C\times\{0\}$. 
But only one of them gives a (TE)-structure on  $M-\KK^{bif}$ 
which extends holomorphically to $M$.
For example, for even $m$ there is a $(TE)$-structure
over $M-\KK^{bif}$ with trivial Stokes structure,
but it does not extend to $M$.

\medskip
(vi) A semisimple $(TEP)$-structure of weight $m\in\Z$ 
over a point $t^0\notin\KK_2$ decomposes formally into 
a sum of rank 1 $(TEP)$-structures of weight $m$,
and for each of them 
the irregular rank 1 factor with $\nnn=d-u_iz^{-2}{\rm d}z$ 
as well as the logarithmic rank 1 factor with 
$\nnn=d+\alpha_i z^{-1}{\rm d}z$ are $(TEP)$-structures
of weight $m$ \cite{HS07}. But the only regular singular rank 1 
$(TEP)$-structure of weight m is the one with
$\alpha_i=m/2$. This forces the values of 
$\AAA^{rse}$ to be all equal to $m/2$. 
This motivates the last part of Conjecture \ref{t7.2} (b). 
\end{remarks}

In section \ref{c8} we will study the 
$(T)$-structures and $(TE)$-structures over the
2-dimensional F-manifolds $I_2(m)$ ($m\geq 3$).
We will prove the following. 

\begin{theorem}\label{t7.4}
Conjecture \ref{t7.2} is true for the germs at 0
of the 2-dimensional F-manifolds of type $I_2(m)$.
\end{theorem}

\begin{remarks}\label{t7.5}
(i) Originally, we hoped to prove Conjecture \ref{t7.2} (b) 
by a constructive and inductive argument looking at
normal forms and their equations,
but up to now, this has been successful only for $I_2(m)$.
It is difficult in this approach to make use of the
generic semisimplicity.

\medskip
(ii) A very different approach is via the theory of 
holonomic $\EE$-modules. 
In \cite{KK81} and \cite{Ka86}, holonomic $\EE$-modules
are considered whose singular support is a smooth 
Lagrangian variety in $T^*M$, and it is shown that they are
{\it simple holonomic $\EE$-modules}.
We hope that this can be used to prove Conjecture \ref{t7.2}
(b) in the case, when the analytic spectrum of 
an irreducible germ of a generically semisimple F-manifold
is smooth. That is precisely the case of the
base space of an isolated hypersurface singularity.
In that case the existence of a 1-parameter family
of $(TE)$-structures is well-known, and the parameter
is the value $\AAA^{rse}(L)\in\C$.
So it remains there only to prove the uniqueness,
and this is the uniqueness of the Stokes structure.

\medskip
(iii) Similarly, we hope that the results in \cite{DDP81}
can be used to prove Conjecture \ref{t7.2} (b) in the case, when
the analytic spectrum of an irreducible germ of a 
generically semisimple F-manifold is a normal
crossing divisor. This would generalize (ii).

\medskip
(iv) Finally, we hope that more elementary arguments
which build on the positive results for the $I_2(m)$
will show at least the uniqueness part in 
Conjecture \ref{t7.2} (b) in the case of an irreducible
germ $(M,t^0)$ of a generically semisimple F-manifold 
with Euler field with the following property: 
At generic points $t^1\in \KK$ 
in the caustic the germ $(M,t^1)$ is isomorphic to 
$A_1^{\dim M-2}I_2(m)$ where $m$ is in a finite subset of 
$\Z_{\geq 3}$. 
One option is to try to control the possible Stokes structures
directly. But also that looks difficult.
\end{remarks}

\section{$(TE)$-structures over the 2-dimensional
F-manifolds $I_2(m)$}\label{c8}
\setcounter{equation}{0}

\noindent 
One aim of this section is to prove Conjecture \ref{t7.2}
in dimension 2, i.e. for the 2-dimensional irreducible
germs of generically semisimple F-manifolds with
Euler fields. Those are the germs of types $I_2(m)$,
$m\in\Z_{\geq 3}$, see Theorem \ref{t2.8} (b).

We will recall in Theorem \ref{t8.2} 
the formal classification of the
(formal or holomorphic) $(T)$-structures over
$I_2(m)$ in \cite[Theorem 16 iii)]{DH19-1}.
Theorem \ref{t8.4} will give semi-normal forms
for holomorphic $(T)$-structures, 
{\it semi} because they are not unique.
It improves \cite[Theorem 16 i)]{DH17}.
It implies Conjecture \ref{t7.2} (a)
in the case of $I_2(m)$. 

Then we turn to $(TE)$-structures over $I_2(m)$
with its Euler field. Theorem \ref{t8.5}
will show that all formal $(TE)$-structures over $I_2(m)$
are formally isomorphic to holomorphic normal forms,
and it will give holomorphic normal forms for all
holomorphic $(TE)$-structures over $I_2(m)$. 
It builds on Theorem \ref{t8.2} (not on Theorem \ref{t8.4})
and uses for the holomorphic classification
Corollary \ref{t5.7}.
In the case $m$ odd, only the well known 1-parameter
family of $(TE)$-structures exists.
In the case $m$ even, a 2-parameter family exists.
It was discovered recently independently in 
\cite{KMS15} and \cite{AL17}.
Finally, Theorem \ref{t8.6} will calculate for
these normal forms the regular singular exponents
(see the Remarks \ref{t7.1}).
Theorem \ref{t8.5}, Theorem \ref{t8.6} and Corollary \ref{t8.7}
together show Conjecture \ref{t7.2} (b) in the case of
$I_2(m)$.

\begin{remarks}\label{t8.1}
(i) We use $t=(t_1,t_2)$ as standard coordinates on $\C^2$,
with coordinate vector fields $\paa_i=\paa/\paa t_i$. 
For a fixed $m\in\Z_{\geq 3}$, 
we shall use the following matrices,
\begin{eqnarray}\label{8.1}
C_1:={\bf 1}_2,\ 
C_2:=\begin{pmatrix}0&t_2^{m-2}\\1&0\end{pmatrix},\ 
D:=\begin{pmatrix}1&0\\0&-1\end{pmatrix},\ 
E:=\begin{pmatrix}0&1\\0&0\end{pmatrix},
\end{eqnarray}
and the relations between them,
\begin{eqnarray}
(C_2)^2=t_2^{m-2}C_1,\ D^2=C_1,\ E^2=0,\label{8.2}\\
C_2D=C_2-2t_2^{m-2}E=-DC_2,\ [C_2,D]=2(C_2-2t_2^{m-2}E),
\label{8.3}\\
C_2E=\frac{1}{2}(C_1-D),\ EC_2=\frac{1}{2}(C_1+D),\ 
[C_2,E]=-D,\label{8.4}\\
DE=E=-ED,\ [D,E]=2E.\label{8.5}
\end{eqnarray}

(ii) The multiplication $\circ$ and the Euler field
$E$ of the F-manifold $I_2(m)$ ($m\in\Z_{\geq 3}$)
on $M=\C^2$ in Theorem \ref{t2.8} (b) are as follows.
Here $\uuuu{\paa}:=(\paa_1,\paa_2)$. 
\begin{eqnarray}\label{8.6}
\paa_1=e,\ \paa_2\circ\paa_2=t_2^{m-2}\paa_1,\quad
\textup{so}\quad
\paa_1\circ\uuuu{\paa}=\uuuu{\paa},\ 
\paa_2\circ\uuuu{\paa}=\uuuu{\paa}\cdot C_2,\\
E=(t_1+c_1)\paa_1 +\frac{2}{m}t_2\paa_2
\textup{ for some }c_1\in\C.\label{8.7}
\end{eqnarray}
We will restrict in the following to $c_1=0$.
This is not a serious restriction, as the F-manifold
structure is constant along the flow of $\paa_1$,
but the Euler field takes up just such a summand 
$c_1\paa_1$ along this flow.

\medskip
(iii) A formal or holomorphic $(T)$-structure or 
$(TE)$-structure over $I_2(m)$ induces via its 
primitive Higgs bundle this multiplication and this
Euler field on $M=\C^2$. That means, that for any
basis $\uuuu{v}$ of it, the matrices 
$A_i=\sum_{k\geq 0}A_i^{(k)}z^k$ and 
$B=\sum_{k\geq 0}B^{(k)}z^k$ with 
\begin{eqnarray}\label{8.8}
\nnn\uuuu{v}=\uuuu{v}\bigl(z^{-1}A_1{\rm d}t_1 + 
z^{-1}A_2{\rm d}t_2+z^{-2}B{\rm d}z\bigr)
\end{eqnarray}
satisfy not only \eqref{4.15} and \eqref{4.16}, but also
\begin{eqnarray}\label{8.9}
A_1^{(0)}=C_1,\ (A_2^{(0)})^2=t_2^{m-2}C_1,\ 
B^{(0)}=-t_1C_1-\frac{2}{m}t_2A_2^{(0)}.
\end{eqnarray}
\end{remarks}

The following theorem from \cite{DH19-1}
gives a unique formal normal form for any formal or 
holomorphic $(T)$-structure
over the germ at 0 of the F-manifold $I_2(m)$.

\begin{theorem}\label{t8.2}
\cite[Theorem 16 iii)]{DH19-1}
Any formal or holomorphic $(T)$-structure 
over the germ at $0$ of the F-manifold $I_2(m)$ 
has a formal basis $\uuuu{v}$ with 
$\nnn\uuuu{v}=\uuuu{v}\cdot (z^{-1}A_1{\rm d}t_1
+z^{-1}A_2{\rm d}t_2)$ and 
\begin{eqnarray}\label{8.10}
&&A_1=C_1,\ A_2=C_2+zf\cdot E,\\
&&\textup{where }
f\ \left\{\begin{array}{ll}
=0&\textup{ if }m=3\\
\in\C[[z]][t_2]_{\leq m-4}&\textup{ if }m\geq 4.
\end{array}\right. \nonumber 
\end{eqnarray}
This $f$ is unique.
\end{theorem}

\begin{remarks}\label{t8.3}
(i) It is not clear whether in the case of a 
holomorphic $(T)$-structure this formal normal
form is holomorphic, so whether then 
$f\in\C\{z\}[t_2]_{\leq m-4}$. More precisely, if we denote 
\begin{eqnarray}
F_1&:=& \C[[z]][t_2]_{\leq m-4},\nonumber\\
F_2&:=& \{f\in F_1\,|\, \textup{the formal }(T)
\textup{-structure for this }f \nonumber \\
&& \textup{is formally isomorphic to a holomorphic }
(T)\textup{-structure}\}\nonumber\\
F_3&:=& \C\{z\}[t_2]_{\leq m-4},\label{8.11}
\end{eqnarray}
we have $F_1\supset F_2\supset F_3$, but we don't know
where $F_2$ is between $F_1$ and $F_3$. 
This is annoying.

(ii) The following theorem gives holomorphic 
semi-normal forms for $(T)$-structures over $I_2(m)$.
{\it Semi}, because they are not unique.

(iii) It is possible to refine them to normal forms by
setting subsets of the coefficients of 
$f\in\C\{t_2\}+ z\cdot\C\{t_2\}$ equal 0. 
The result is combinatorially surprisingly rich, 
but still not appealing. 
\end{remarks}

The proof of Theorem \ref{t8.4}
follows roughly the proof of \cite[Theorem 16 i)]{DH19-1},
but it uses also Theorem \ref{t5.1} (a). 
Therefore it is simpler and gives stronger semi-normal forms.

\begin{theorem}\label{t8.4} 
Over the F-manifold $I_2(m)$ for $m\geq 3$, 
any holomorphic (T)-structure 
has a holomorphic basis $\uuuu{v}$ with 
$\nnn\uuuu{v}=\uuuu{v}\cdot (z^{-1}A_1{\rm d}t_1
+z^{-1}A_2{\rm d}t_2)$ and 
\begin{eqnarray}\label{8.12}
A_1 &=& A_1^{(0)}=C_1,\\
A_2 &=& \sum_{k=0}^2 A_2^{(k)} = C_2+zf\cdot E
\label{8.13}\\
&&\textup{with }f\in \C\{t_2\} + z\cdot\C\{t_2\}.\nonumber
\end{eqnarray}
\end{theorem}

{\bf Proof:} 
Extend the (T)-structure to a pure (TL)-structure,
using Theorem \ref{t5.1}, 
and follow the steps 1 to 3 in the proof of 
\cite[Theorem 16 i)]{DH19-1}. 
The details are as follows.

Choose a basis $\uuuu{v}=(v_1,v_2)$ of global sections
whose restrictions to $\{\infty\}\times M$ are $\nabla^{res,\infty}$-flat 
and with $v_1|_{\{0\}\times M}$ primitive. Then 
\begin{eqnarray}\label{8.14}
\nabla \uuuu{v}&=&\uuuu{v}\cdot
(z^{-1}A_1^{(0)}{\rm d}t_1+z^{-1}A_2^{(0)}{\rm d}t_2)\\
\textup{where } A_1^{(0)}&=&C_1,
\quad (A_2^{(0)})^2=t_2^{m-2}C_1,\nonumber
\end{eqnarray}
and a priori $A_2^{(0)}\in M_{2\times 2}(\C\{t_1,t_2\})$,
but \eqref{5.3}, 
$\paa_1 A_2^{(0)}=\paa_2A_1^{(0)}=\paa_2 C_1=0$, 
shows immediately $A_2^{(0)}\in M_{2\times 2}(\C\{t_2\})$.

There are unique vector fields $\www\paa_1$ and 
$\www\paa_2$ with 
$C_{\www\paa_j}v_1|_{\{0\}\times M}=v_j|_{\{0\}\times M}$
and thus $z\nabla_{\www\paa_j}v_1 =v_j$. 
Then
\begin{eqnarray}\label{8.15}
\paa_2\circ (\www \paa_1,\www\paa_1)&=&
(\www\paa_1,\www\paa_2)\cdot A_2^{(0)}.
\end{eqnarray}
Of course $\www\paa_1=\paa_1$.
We can assume that $v_2$ is chosen so that 
$\www\paa_2|_{t=0}=\paa_2|_{t=0}$.
Then for suitable $a,b\in\C\{t_2\}$ with $b(0)=1$ and $a(0)=0$ 
\begin{eqnarray}\label{8.16}
\www\paa_2 &=& ab^{-1}\paa_1 + b^{-1}\paa_2,\\
\textup{so }(\www\paa_1,\www\paa_2)&=& (\paa_1,\paa_2)\cdot 
\begin{pmatrix}1& ab^{-1}\\ 0 &b^{-1}\end{pmatrix},
\nonumber\\
(\paa_1,\paa_2) &=& (\www\paa_1,\www\paa_2)\cdot T\quad 
\textup{with }T:= \begin{pmatrix}1&-a\\0&b\end{pmatrix}
=\begin{pmatrix}1& ab^{-1}\\ 0 &b^{-1}\end{pmatrix}^{-1}.
\nonumber
\end{eqnarray}
Choose the new basis
\begin{eqnarray}\label{8.17}
\uuuu{\www v}=(\www v_1,\www v_2)
&:=& \uuuu{v}\cdot T
\quad \textup{with }T:= \begin{pmatrix}1&-a\\0&b\end{pmatrix}.
\end{eqnarray}
Its restriction to $\{\infty\}\times M$ is no longer
$\nnn^{res,\infty}$-flat. But the matrices $\www A_1$
and $\www A_2$ 
with $\nnn\uuuu{\www v}=\uuuu{\www v}\cdot
(z^{-1}\www A_1{\rm d}t_1+z^{-1}\www A_2{\rm d}t_2)$ satisfy
$\www A_1=C_1$ and 
\begin{eqnarray}
\www A_2 &=& T^{-1}z\paa_2 T + T^{-1}A_2 T \nonumber\\
&=& ... = \begin{pmatrix}0& t_2^{m-2}\\ 1& 0\end{pmatrix}
+ z\cdot\begin{pmatrix}0 & -\paa_2 a +ab^{-1}\paa_2 b\\ 0 & b^{-1}\paa_2 b
\end{pmatrix}\label{8.18}\\
&=& C_2 + z\cdot b^{-1}\paa_2 b\cdot\frac{1}{2}(C_1-D)
+ z\cdot (-\paa_2 a + ab^{-1}\paa_2b)\cdot E.\nonumber
\end{eqnarray}

By step 2 in the proof of Theorem 16 i) in \cite{DH19-1}, we
can change $\uuuu{\www v}$ to another basis so that its matrix 
$A_2^{new,1}$ has the form
\begin{eqnarray}\label{8.19}
A_2^{new,1}
&=&C_2 + z\cdot b^{-1}\paa_2 b\cdot\frac{-1}{2}\cdot D
 + z\cdot (-\paa_2 a + ab^{-1}\paa_2b)\cdot E\hspace*{1cm}\\
 &=:& C_2 + z\cdot a_3\cdot D + z\cdot a_4\cdot E.\nonumber
\end{eqnarray}

Now one can apply step 3 in the proof of Theorem 16 i) in
\cite{DH19-1} with $\tau_3:=0$ and $\tau_4:=a_3$, so with 
\begin{eqnarray}\label{8.20}
\www T&:=& C_1+za_3 E.
\end{eqnarray}
The new matrix $A_2^{new,2}$ has the shape
\begin{eqnarray}
A_2^{new,2}&=& {\www T}^{-1}z\paa_2{\www T} 
+ {\www T}^{-1}A_2^{new,1}{\www T}
=C_2+zf\cdot E \nonumber \\
\textup{with }f&=& a_4+z(\paa_2a_3+a_3^2)
\in\C\{t_2\} + z\cdot\C\{t_2\}.
\label{8.21}
\end{eqnarray}
(In the notation of \cite{DH19-1} $\www a_1=0,f=\www a_4$;
details of the calculation leading to \eqref{8.21}
are given in \cite{DH19-1}.)
\hfill $\Box$

\bigskip

Now we come to $(TE)$-structures over the F-manifold
$I_2(m)$ with the Euler field in \eqref{8.7} with $c_1=0$. 
Theorem \ref{t8.5} will show that any formal $(TE)$-structure
is formally isomorphic to a holomorphic $(TE)$-structures,
and it will give unique holomorphic normal forms for the
holomorphic $(TE)$-structures. They have 2 parameters 
$(\alpha,\lambda)\in\C^2$ if $m$ is even and 1 parameter 
$\alpha\in\C$ if $m$ is odd.

\begin{theorem}\label{t8.5}
Any formal $(TE)$-structure over the germ at $0$ of the
F-manifold $I_2(m)$ with the Euler field 
$E=t_1\paa_1+\frac{2}{m}t_2\paa_2$ is formally isomorphic to
a holomorphic $(TE)$-structure. 
Any holomorphic $(TE)$-structure has a basis $\uuuu{v}$
such that the matrices $A_1,A_2$ and $B$ in \eqref{8.8}
are in the following normal form.
\begin{eqnarray}\label{8.22}
A_1&=& C_1,\\
A_2&=& C_2 + z\cdot f\cdot E,\label{8.23}\\
B &=& \bigl(-t_1C_1-\frac{2}{m}t_2C_2\bigr) \label{8.24} \\
&+& z\cdot\Bigl(\alpha\cdot C_1+\frac{2-m}{2m}\cdot D
+\frac{-2}{m}t_2f\cdot E\Bigr),\nonumber\\
\textup{with }
f&=& \left\{\begin{array}{ll}
0 &\textup{if }m\textup{ is odd,}\\
\lambda\cdot t_2^{(m-4)/2}&\textup{if }m\textup{ is even.}
\end{array}\right. \label{8.25}
\end{eqnarray}
Here $\alpha\in\C$ and $\lambda\in\C$. Any formal or
holomorphic $(TE)$-structure has a unique such normal form. 
\end{theorem}

{\bf Proof:}
Let a formal or holomorphic $(TE)$-structure over the 
$F$-manifold $I_2(m)$ be given. Because of Theorem \ref{t8.2}, 
we can choose a formal basis $\uuuu{v}$ with
the matrices $A_1$ and $A_2$ as in \eqref{8.10}, so 
\begin{eqnarray}\label{8.26}
A_1=C_1\textup{ and }A_2=C_2+zf\cdot E
\end{eqnarray}
with $f\in\C[[z]][t_2]_{\leq m-4}$ if $ m\geq 4$ 
and $f=0$ if $m=3$.

The main part of this proof will show that this basis
can be changed to a formal basis with matrices 
$A_1,A_2$ and $B$ as in Theorem \ref{t8.5}.

At the end of the proof, we can and will apply 
Corollary \ref{t5.7}.
That will show that a holomorphic $(TE)$-structure
has a holomorphic basis with the normal form in
Theorem \ref{t8.5}.

We will show that the matrix $A_2$ satisfies 
\eqref{8.23} and \eqref{8.25}
and that the matrix $B$ can be brought by a base change,
which does not change $A_1$ and $A_2$, to the form in 
\eqref{8.24}.

Equation \eqref{4.16} for $i=1$ gives
\begin{eqnarray}\label{8.27}
0=z\paa_1 B-z^2\paa_z A_1+[A_1,B]=z\paa_1 B.
\end{eqnarray}
Because of \eqref{8.9}, $B^{(0)}=-t_1C_1-\frac{2}{m}C_2$. 
Therefore $B$ has the form
\begin{eqnarray}\label{8.28}
B&=& (-t_1+zb_1)C_1+b_2C_2+zb_3D+zb_4E\\
&&\textup{with }b_1,b_2,b_3,b_4\in\C\{t_2,z]].\nonumber
\end{eqnarray}
We know already $b_2\in -\frac{2}{m}t_2+z\C\{t_2,z]]$,
but we will recover this also below. 
Equation \eqref{4.16} for $i=2$ gives strong constraints
on $b_1,b_2,b_3$ and $b_4$. 
\begin{eqnarray*}
0&=& z\paa_2 B-z^2\paa_zA_2+zA_2+[A_2,B]\\
&=& z^2\paa_2b_1\cdot C_1+z\paa_2b_2\cdot C_2+zb_2\cdot(m-2)t_2^{m-3}\cdot E\\
&&+z^2\paa_2b_3\cdot D+z^2\paa_2b_4\cdot E\\
&+& z\cdot C_2+(-z^2\paa_z(z\cdot f)+z^2f)E\\
&+& z\cdot b_3\cdot 2(C_2-2t_2^{m-2}E)+z\cdot b_4\cdot (-D)\\
&+& z\cdot fb_2\cdot D+z^2\cdot fb_3\cdot (-2E)\\
&=& z^2\cdot C_1\cdot (\paa_2b_1)\\
&+&z\cdot C_2\cdot (\paa_2b_2+1+2b_3)\\
&+& z\cdot D\cdot (-b_4+fb_2+z\cdot \paa_2b_3)\\
&+& z\cdot E\cdot ((m-2)t_2^{m-3}\cdot b_2 + z\cdot \paa_2b_4\\
&&-z^2\paa_z f-4t_2^{m-2}b_3-z\cdot 2fb_3).
\end{eqnarray*}
The coefficients of $C_1,C_2,D$ and $E$ give
\begin{eqnarray}\label{8.29}
\paa_2b_1&=&0,\quad\textup{i.e. }b_1\in \C[[z]],\\
2b_3&=& -\paa_2b_2-1,\label{8.30}\\
b_4 &=& fb_2+z\cdot \paa_2b_3 = fb_2-z\cdot \frac{1}{2}\paa_2^2b_2,\label{8.31}\\
0&=& t_2^{m-3}((m-2)\id +2t_2\paa_2)(b_2)+2t_2^{m-2}
\label{8.32}\\
&+& z\cdot(-z\paa_z f+\paa_2f\cdot b_2+2f\cdot\paa_2b_2+f
-z\cdot \frac{1}{2}\paa_2^3 b_2).\nonumber
\end{eqnarray}
The following choice of $T$ reduces $b_1\in \C[[z]]$ to 
$\www b_1=z\cdot\alpha$ for some $\alpha\in\C$ 
without changing anything else. Set
\begin{eqnarray}
b_1&=& \alpha+\sum_{k\geq 1}b_1^{(1)}\cdot z^k,
\nonumber\\
T&:=&\exp\left(-\sum_{k\geq 1}\frac{b_1^{(k)}}{k}
\cdot z^{k}\right)\cdot C_1, \label{8.33}\\
\textup{thus }z\paa_z T&=& (\alpha-b_1)C_1\cdot T.\nonumber
\end{eqnarray}
\eqref{4.25} and \eqref{4.26} show that $A_1$ and $A_2$ 
are unchanged and that the only change in $B$ is that $b_1$ 
is replaced by $\www b_1=\alpha\cdot z$.

Write 
$$b_2=\sum_{k\geq 0}b_2^{(k)}(t_2)\cdot z^k
\quad\textup{and}\quad
f=\sum_{k\geq 0}f^{(k)}(t_2)\cdot z^k.$$ 
Then $b_2^{(0)}$ is determined uniquely by the first line 
of \eqref{8.32}, and it is 
$$b_2^{(0)}=\frac{-2}{m}t_2.$$ 
The next term $b_2^{(1)}$ is determined by both lines 
of \eqref{8.32}, but in the second line only $b_2^{(0)}$ and 
$f^{(0)}$ are relevant. 
Because all $f^{(k)}\in\C[t_2]_{\leq m-4}$
and $b_2^{(0)}\in\C[t_2]_{\leq 1}$, the first line 
for $b_2^{(1)}$ and the second line for $b_2^{(0)}$ 
and $f^{(0)}$ must vanish separately. This gives 
$$b_2^{(1)}=0.$$
By the same argument, one obtains inductively
\begin{eqnarray}\label{8.34}
b_2^{(k)}=0\quad\textup{for any }k\geq 1,\quad\textup{thus }
b_2=b_2^{(0)}=\frac{-2}{m}t_2
\end{eqnarray}
Now the second line of \eqref{8.32} becomes
\begin{eqnarray*}
0&=& (mz\paa_z + 2t_2\paa_2 + (4-m)\id)(f).
\end{eqnarray*}
The only solutions are
\begin{eqnarray}\label{8.35}
f=f^{(0)} =\left\{\begin{array}{ll}
\lambda\cdot t_2^{(m-4)/2}\hspace*{1cm}&
\textup{if }m\textup{ is even,}\\
0 &\textup{if }m\textup{ is odd.}\end{array}\right. 
\end{eqnarray}
for an arbitrary $\lambda\in\C$ if $m$ is even. Then
\begin{eqnarray}\label{8.36}
b_3 &=& -\frac{1}{2}\paa_2b_2-\frac{1}{2}=\frac{2-m}{2m},\\
b_4 &=& fb_2=\frac{-2}{m}t_1\cdot f^{(0)}
=\frac{-2}{m}\cdot\lambda\cdot t_2^{(m-2)/2}.
\label{8.37}
\end{eqnarray}
This gives \eqref{8.23}, \eqref{8.24} and \eqref{8.25}.

In the case of a formal $(TE)$-structure we are ready.
In the case of a holomorphic $(TE)$-structure,
one observes that the matrix $B$ restricted to $t=0$
gives a $(TE)$-structure over $t^0=0$ with a regular
singular pole (even a logarithmic pole) at $z=0$. 
By Corollary \ref{t5.7}, the base change matrix 
$\www T$ from any holomorphic basis of the original
$(TE)$-structure to the formal basis $\uuuu{v}$
above is holomorphic.
Therefore the basis $\uuuu{v}$ is a holomorphic basis
of the $(TE)$-structure. 

The uniqueness of the normal form in Theorem \ref{t8.5} 
follows from the uniqueness of the normal form in 
Theorem \ref{t8.2} (the proof did not change $f$) 
and from the uniqueness of $\alpha$ above.
\hfill$\Box$

\bigskip

For the proof of Conjecture \ref{t7.2} (b), 
we need the regular singular exponents 
of the normal forms in Theorem \ref{t8.5}.
Theorem \ref{t8.6} gives these values
and shows the bijection between the normal forms
and the possible maps 
$\AAA:\{\textup{components of }L_M\}\to\C$
from \eqref{7.2}.
Together with the uniqueness of the normal form in
Theorem \ref{t8.5}, this shows Conjecture \ref{t7.2} (b)
for the $I_2(m)$, up to the statement on
$(TEP)$-structures. That follows from Corollary \ref{t8.7}.

\begin{theorem}\label{t8.6}
Consider a $(TE)$-structure over the germ at 0 of the
F-manifold $I_2(m)$ ($m\in\Z_{\geq 3}$).
Let $\uuuu{v}$ be a basis with matrices $A_1,A_2,B$
in the normal form in Theorem \ref{t8.5}.

(i) Suppose that $m$ is even. Then $L_M$ has two components.
And then $\alpha+\frac{\mp\lambda}{m}$ 
are the two regular singular exponents. 

(ii) Suppose that $m$ is odd. Then $L_M$ has one component.
And then $\alpha$ is the only regular singular exponent.
\end{theorem}

{\bf Proof:} 
Locally on $M-\KK$ choose for odd $m$ 
a square root $t_2^{1/2}$ of $t_2$, 
define (for even and for odd $m$) 
$\tau:=t_2^{(m-2)/2}$, and consider the base change matrices
\begin{eqnarray}\label{8.38}
T=T^{(0)}= \begin{pmatrix} \tau & -\tau\\ 1 & 1 \end{pmatrix}
\textup{ and }
T^{-1}=\frac{1}{2\tau}\begin{pmatrix}1&\tau\\ -1& \tau \end{pmatrix}.
\end{eqnarray} 
Conjugation with $T$ diagonalizes $C_2$
and transforms $D$ and $E$ as follows,  
\begin{eqnarray}\label{8.39}
T^{-1}C_2T &=& \begin{pmatrix}\tau&0\\ 0& -\tau \end{pmatrix},\\
T^{-1}DT &=& \begin{pmatrix}0&-1\\-1&0 \end{pmatrix},
\label{8.40}\\
T^{-1}ET &=& \frac{1}{2\tau}
\begin{pmatrix}1&1\\-1& -1 \end{pmatrix},\label{8.41}\\
\textup{and }T^{-1}\paa_2 T
&=& \frac{m-2}{4t_2}\cdot\begin{pmatrix}1&-1\\-1& 1
\end{pmatrix}. \label{8.42}
\end{eqnarray}

The local basis
\begin{eqnarray}\label{8.43}
\uuuu{\www v}:=\uuuu{v}\cdot T
\end{eqnarray}
satisfies 
\begin{eqnarray}
\nnn \uuuu{\www v}&=& \uuuu{\www v}
\cdot(z^{-1}\www A_1\ddd t_1 +z^{-1}\www A_2\ddd t_2+
z^{-2}\www B\ddd z)
\qquad\textup{with}\nonumber\\
\www A_1 &=& C_1,\nonumber\\
\www A_2 &=& T^{-1}z\paa_2 T+T^{-1}A_2T\\
&=& \begin{pmatrix}\tau&0\\ 0& -\tau \end{pmatrix}\nonumber 
+z\frac{m-2}{4t_2}\cdot\begin{pmatrix}1&-1\\-1& 1 \end{pmatrix}
+z\frac{f}{2\tau}\cdot\begin{pmatrix}1&1\\-1& -1 \end{pmatrix}
\hspace*{1cm}\label{8.44}
\end{eqnarray}
\begin{eqnarray}
\www B&=& T^{-1}z^2\paa_z T+T^{-1}BT=T^{-1}BT \nonumber\\
&=& \left(-t_1C_1-\frac{2}{m}t_2\tau D\right) 
+ z\cdot\begin{pmatrix}\alpha -\frac{\lambda}{m}&-\frac{2-m}{2m} -\frac{\lambda}{m}\\
-\frac{2-m}{2m} +\frac{\lambda}{m}&\alpha +\frac{\lambda}{m}\end{pmatrix}. \hspace*{1cm}\label{8.45}
\end{eqnarray}
The matrix ${\www B}^{(0)}$ is diagonal.
Therefore the diagonal entries of the matrix 
${\www B}^{(1)}$ are the regular singular exponents.
So, they are $\alpha+\frac{\mp\lambda}{m}$. 
This proves Theorem \ref{t8.6} in both cases, 
whether $m$ is even or odd.
\hfill$\Box$

\begin{corollary}\label{t8.7}
Consider a $(TE)$-structure over the germ at 0 of the
F-manifold $I_2(m)$ ($m\in\Z_{\geq 3}$).
Let $\uuuu{v}$ be a basis with matrices $A_1,A_2,B$
in the normal form in Theorem \ref{t8.5}.
Let $w\in\Z$. 

The $(TE)$-structure extends to a $(TEP)$-structure
of weight $w\in\Z$ if and only if $\lambda=0$ (if $m$ is even)
and $\alpha=\frac{w}{2}$. The pairing $P$ is unique
up to rescaling and is given by
\begin{eqnarray}\label{8.46}
z^{-w}\cdot P(\uuuu{v}^t,\uuuu{v})=P^{mat,(0)}=
c\cdot\begin{pmatrix}0&1\\1&0\end{pmatrix}
\end{eqnarray}
with $c\in\C^*$ arbitrary. 
\end{corollary}

{\bf Proof:} Remark \ref{t7.3} (vi) and Theorem \ref{t8.6}
together show that $\alpha+\frac{\mp\lambda}{m}=\frac{w}{2}$
is necessary for the extendability of the $(TE)$-structure
to a $(TEP)$-structure of weight $w$. Consider
the case $\lambda=0$ (if $m$ is even) and $\alpha=\frac{w}{2}$.
If a pairing $P$ exists, it is given by a matrix
$P^{mat,(0)}=z^{-w}P(\uuuu{v}^t,\uuuu{v})$ which satisfies
\eqref{5.8}--\eqref{5.11} and \eqref{4.14}.
Thus it is a constant nondegenerate and symmetric matrix
with 
\begin{eqnarray}\label{8.47}
A_2^{(0)} P^{mat,(0)}-P^{mat,(0)}A_2^{(0)}=0,\\
(-w)P^{mat,(0)}+ (B^{(1)})^t P^{mat,(0)} 
+ P^{mat,(0)}B^{(1)}=0.
\label{8.48}
\end{eqnarray}
and $A_2^{(0)}=C_2$, $B^{(1)}=\frac{w}{2}C_1+\frac{2-m}{2m}D$.
\eqref{8.47} alone forces $P^{mat,(0)}$ 
to be as in \eqref{8.46}.
\eqref{8.48} is then satisfied. \hfill$\Box$

\begin{remark}\label{t8.8}
In Remark \ref{t3.4} (iii) it was claimed that one
has above the F-manifold $I_2(m)$ with Euler field
only a 2-dimensional (if $m$ is even) respectively
1-dimensional (if $m$ is odd) family of flat F-manifolds with
Euler field. This is close to Theorem \ref{t8.5},
but stronger. We still have to make the steps (III)(b)
and (IV)(b) (in the notation of Remark \ref{t6.11}).
And we have to see that no new parameters arise 
during these steps.

The step (III) (b): 
The restriction of the $(TE)$-structure in 
Theorem \ref{t8.5} to $t^0=0$ has a logarithmic pole
at $z=0$ with residue eigenvalues 
$\alpha\pm \frac{2-m}{2m}$. Their difference is smaller
than 1. Therefore this $(TE)$-structure over $t^0=0$
has a unique extension to a pure $(TLE)$-structure.
And therefore also the whole $(TE)$-structure has a 
unique extension to a pure $(TLE)$-structure (Theorem \ref{t5.1} (c)).
In the step (III) (b) no new parameter arises.
If $m$ is odd or if $m$ is even and $\lambda=0$, 
the basis $\uuuu{v}$ in Theorem \ref{t8.5} 
is already a basis of the pure $(TLE)$-structure.
If $m$ is even and $\lambda\neq 0$, 
the following base change leads to 
the pure $(TLE)$-structure. Define 
\begin{eqnarray}\label{8.49}
\beta&:=& \frac{2}{2-m}\lambda t_2^{(m-2)/2},\\
T&:=& T^{(0)}=\begin{pmatrix}1&\beta\\0&1\end{pmatrix}\label{8.50},\\
\uuuu{\www v}&:=& \uuuu{v}\cdot T.\label{8.51}
\end{eqnarray}
Then
\begin{eqnarray}\label{8.52}
\www A_1&=& T^{-1}A_1T+T^{-1}z\paa_1 T = C_1,\\
\www A_2&=& T^{-1}A_1T+T^{-1}z\paa_2 T =
\begin{pmatrix}-\beta & t_2^{m-2}-\beta^2\\1&\beta\end{pmatrix},
\label{8.53}\\
\www B&=& T^{-1}BT+T^{-1}z^2\paa_z T \nonumber\\
&=& -t_1C_1-\frac{2}{m}t_2\www A_2 
+ z\bigl(\alpha\cdot C_1 + \frac{2-m}{2m}\cdot D\bigr).
\label{8.54}
\end{eqnarray}

The step (IV) (b): 
Theorem \ref{t6.6} (b) shows how to make this step.
One needs a section $\omega\in\C\cdot {\www v_1}+
\C\cdot {\www v_2}$ such that its restriction
$\zeta:=\omega|_{\{0\}\times M}$ to $K$ is a primitive
section of the Higgs bundle and an eigenvector of $Q$.
Here $Q$ is given by the matrix $-{\www B}^{(1)}$,
which is diagonal. Thus the eigenvectors of $Q$ are
${\www v_1}$ and ${\www v_2}$. 
The second one, ${\www v_2}$ is not a primitive section,
as the second column of the matrix ${\www A}_2={\www A}_2^{(0)}$ 
vanishes at $t^0=0$. Only the first one, ${\www v_1}$
is a primitive section. Also in the step (IV) (b) 
no new parameter arises. 

The flat coordinates $\www t_1,\www t_2$ and the flat vector fields
$\www\paa_1,\www\paa_2$ can be made explicit easily.
$C_{\www\paa_i}(\www v_1|_{\{0\}\times M})
=\www v_i|_{\{0\}\times M}$ leads to 
\begin{eqnarray}\label{8.55}
(\www \paa_1,\www \paa_2)&=& (\paa_1,\paa_2+\beta\paa_1),\\
(\www t_1,\www t_2)&=& (t_1-\frac{2}{m}t_2\beta,t_2),\label{8.56}\\
\www\paa_2\circ\www\paa_2 &=& (t_2^{m-2}+\beta^2)\paa_1 + 2\beta\paa_2
= (t_2^{m-2}-\beta^2)\www\paa_1 + 2\beta\www\paa_2.
\hspace*{0.5cm}\label{8.57}
\end{eqnarray}
A vector potential (see Remark \ref{t3.2} (iv)) is
\begin{eqnarray}\label{8.58}
\Bigl(\frac{1}{2}{\www t_1}^2 + 
\frac{1}{m(m-1)}(t_2^m-t_2^2\beta^2)\Bigr)\www\paa_1 
+ \Bigl(\www t_1t_2+\frac{8}{(m+2)m}t_2^2\beta\Bigr)\www\paa_2.
\end{eqnarray}
\end{remark}

\end{document}